\documentclass[journal]{IEEEtran}%

\usepackage{theorem}%
\usepackage{amsmath,amssymb,amsfonts}%
\usepackage{graphicx}      % include this line if your document contains figures
\usepackage{color,dsfont}
\usepackage{enumerate,upgreek}
\definecolor{mycol}{rgb}{0,0,1}
\definecolor{mcc}{rgb}{0,0.4,0.6}
\usepackage{amsmath,amssymb}
\usepackage[mathscr]{eucal}
\usepackage{times,pifont}
\usepackage{framed}
\usepackage{enumitem}
\setlist[itemize]{leftmargin=*}
\usepackage[bookmarks,bookmarksopen,bookmarksdepth=3]{hyperref}
\usepackage{xcolor}

\usepackage{csquotes}  % For \enquotes
\newcommand\q{\enquote}

\usepackage{tikz}
\usetikzlibrary{arrows}
\usetikzlibrary{graphs,decorations.pathmorphing,decorations.markings}
\usetikzlibrary{calc}

\newtheorem{definition}{Definition}[section]%
\newtheorem{proposition}[definition]{Proposition}%
\newtheorem{theorem}[definition]{Theorem}%
\newtheorem{lemma}[definition]{Lemma}%
\newtheorem{assumption}[definition]{Assumption}%

{\theorembodyfont{\rmfamily}\newtheorem{remark}[definition]{Remark}}%
{\theorembodyfont{\rmfamily}\newtheorem{example}[definition]{Example}}%

\newenvironment{proof}
  {{\bf Proof:}}
  {\qquad \hspace*{\fill} $\Box$}%

\newcommand{\unit}{\mathds{1}}%

\newcommand{\N}{\mathbb{N}}%
\newcommand{\R}{\mathbb{R}}%
\newcommand{\Z}{\mathbb{Z}}%
\newcommand{\id}{\mathrm{id}}%
\newcommand{\KC}{\mathcal{K}}%
\newcommand{\LC}{\mathcal{L}}%
\newcommand{\UC}{\mathcal{U}}%
\newcommand{\VC}{\mathcal{V}}%

\newcommand{\PC}{\mathcal{P}}%
\newcommand{\tm}{\times}%
\newcommand{\rmD}{\mathrm{D}}%
\newcommand{\rmd}{\mathrm{d}}%
\newcommand{\inner}{\mathrm{int}}%
\newcommand{\ep}{\varepsilon}%
\newcommand{\esssup}{\operatorname*{ess\;sup}}%
\newcommand{\K}{\mathcal{K}}%
\newcommand{\Kinf}{\mathcal{K_\infty}}%
\newcommand{\KL}{\mathcal{KL}}%
\newcommand{\LL}{\mathcal{L}}%

\newcommand{\im}{\mathrm{im}}%

\allowdisplaybreaks%

\begin{document}

\title{A Lyapunov-based ISS small-gain theorem for infinite networks of nonlinear systems}%
%\author{Christoph Kawan, Andrii Mironchenko, Navid Noroozi, Majid Zamani
\author{Christoph Kawan, Andrii Mironchenko, Majid Zamani,~\IEEEmembership{Senior Member,~IEEE}

	\thanks{C.~Kawan is with the Institute of Informatics, LMU Munich, Germany; e-mail: \texttt{christoph.kawan@lmu.de}. His work is supported by the German Research Foundation (DFG) through the grant ZA 873/4-1.}%
  \thanks{A.~Mironchenko is with Faculty of Computer Science and Mathematics, University of Passau, 94032 Passau, Germany; e-mail: \texttt{andrii.mironchenko@uni-passau.de}. His work is supported by the DFG through the grant MI 1886/2-1.}
	\thanks{M.~Zamani is with the Computer Science Department, University of Colorado Boulder, CO 80309, USA. 
M.~Zamani is also with the the Institute of Informatics, LMU Munich, Germany; email: {\tt\small majid.zamani@colorado.edu}. His work is supported in part by the DFG through the grant ZA 873/4-1 and the H2020 ERC Starting Grant AutoCPS (grant agreement No.~804639).}
}
\date{}%
\maketitle%

\begin{abstract}
In this paper, we show that an infinite network of input-to-state stable (ISS) subsystems, admitting ISS Lyapunov functions, itself admits an ISS Lyapunov function, provided that the couplings between the subsystems are sufficiently weak. The strength of the couplings is described in terms of the properties of an infinite-dimensional nonlinear positive operator, built from the interconnection gains. If this operator induces a uniformly globally asymptotically stable (UGAS) system, a Lyapunov function for the infinite network can be constructed. We analyze necessary and sufficient conditions for UGAS and relate them to small-gain conditions used in the stability analysis of finite networks. 
\end{abstract}

\begin{IEEEkeywords}
Nonlinear systems, small-gain theorems, infinite-dimensional systems, input-to-state stability, Lyapunov methods, large-scale systems%
\end{IEEEkeywords}

\section{Introduction}

We are surrounded by networks: social networks, power grids, transportation and manufacturing networks, etc. These networks grow in size every year, and emerging technologies, such as Cloud Computing and 5G communication, make this trend even more apparent. As stability properties of networks may deteriorate with the increase in the number of participating agents \cite{Bamieh.2012}, it is natural to study infinite networks which over-approximate large-scale networks as a worst-case scenario.%

A prominent place in these investigations is occupied by the theory of linear \emph{spatially invariant systems} \cite{BPD02,BaV05,CIZ09}. In such networks, infinitely many subsystems are coupled via the same pattern. This nice geometrical structure together with the \emph{linearity} of subsystems allowed researchers to develop powerful criteria for the stability of such networks.%

On the other hand, in the stability analysis of \emph{finite networks with nonlinear components}, groundbreaking results have been obtained within the framework of \emph{input-to-state stability (ISS)} \cite{Son89}. According to the ISS \emph{small-gain approach}, the influence of any subsystem on other subsystems in a network is characterized by so-called gain functions. The \emph{gain operator} constructed from these functions characterizes the interconnection structure of the network. The small-gain theorems for finite networks of ISS systems described by ordinary differential equations (ODEs) \cite{JTP94,JMW96,DRW07,Dashkovskiy.2010} show that if the gains are small enough (which is expressed in terms of a so-called \emph{small-gain condition} of the gain operator), the network is ISS. These results have numerous applications in systems theory \cite{LJH14,DES11}, and play a major role in a large part of modern nonlinear control \cite{Son08,KKK95}.%

Recently, significant advances have been achieved in an \emph{infinite-dimensional ISS theory}, see \cite{KaK19,MiP20} for a comprehensive overview of the topic, and \cite{Sch20} for an overview of the linear theory.%

This progress motivated the development of the ISS small-gain framework for the stability analysis of infinite interconnections of nonlinear systems without any spatial invariance assumption. This research was initiated in \cite{DaP19}, where nonlinear Lyapunov-based small-gain theorems have been obtained under the very strong assumption that all gains are uniformly less than the identity. In \cite{DaP19}, the authors also apply their small-gain theorem for stabilization of infinite nonlinear networks by revisiting a backstepping method. In \cite{KMS21}, tight Lyapunov-based small-gain theorems have been obtained for networks of exponentially ISS systems with linear gains, and these results have been applied to distributed observer design and cooperative control of infinite networks in \cite{NMK20b}.%

Nonlinear \emph{trajectory-based small-gain theorems} for infinite networks have been developed in \cite{MKG20}, where it was shown that an infinite network of ISS systems is ISS if the corresponding nonlinear gain operator satisfies the so-called \emph{monotone limit property}. The monotone limit property implies the \emph{uniform small-gain condition} \cite{MKG20}, which is equivalent to the monotone bounded invertibility property. The latter played a key role in the derivation of the ISS small-gain theorem for finite networks in \cite[Lem.~13]{DRW07}.%

This paper is strongly motivated by \cite{dashkovskiy2019stability}, where the \emph{robust strong small-gain condition} was introduced and a method to construct paths of strict decay, based on the concept of the \emph{strong transitive closure} of the gain operator, was proposed. For finite networks, this method was used in \cite[Prop.~2.7, Rem.~2.8]{KaJ11}, see also \cite{Rue17} for more details on the importance of this concept in small-gain theory. Based on these results, in \cite{dashkovskiy2019stability} a small-gain theorem for infinite networks, and the construction of associated ISS Lyapunov functions were proposed under the assumption that a \emph{linear} path of strict decay for the gain operator exists. Although there are examples of nonlinear infinite networks with nonlinear gains, whose gain operators admit a linear path of strict decay, in general, this requirement is quite restrictive. In fact, the Lyapunov-based small-gain theorems for finite networks developed in \cite{Dashkovskiy.2010} do not require the linearity of the path of strict decay.%

\textbf{Contribution.} We consider infinite networks of ISS control systems described by ODEs, admitting ISS Lyapunov functions with a corresponding gain operator, characterizing the influence of the subsystems on each other. In our first main result (Theorem~\ref{thm_mainres}), we show that \emph{the existence of a (possibly nonlinear) path of strict decay for the gain operator $\Gamma$ (together with some uniformity conditions) implies that the whole network is ISS and a corresponding ISS Lyapunov function for the network can be constructed}. Our result extends (up to some extra uniformity condition which we require from the path of strict decay) the nonlinear Lyapunov-based small-gain theorem for finite networks (in maximum formulation) shown in \cite{Dashkovskiy.2010} to the setting of infinite networks, recovers the Lyapunov-based small-gain theorem for infinite networks from \cite{dashkovskiy2019stability}, and partially recovers the main result in \cite{DaP19}; see Section~\ref{sec:Discussion of obtained results} for a detailed discussion.%

Next, we introduce the concept of a \emph{max-robust small-gain condition}, which is less conservative than the robust small-gain condition from \cite{dashkovskiy2019stability}, but better compatible with max-type gain operators, and can be fully characterized in terms of the asymptotic properties of the discrete-time system induced by the gain operator.%

In our second main result (Theorem~\ref{thm_omegapath}), we show that the \emph{uniform global asymptotic stability (UGAS) of the system induced by a scaled gain operator guarantees the existence of a path of strict decay}. Furthermore, we explicitly construct this path via the concept of the strong transitive closure of the gain operator, adopted from \cite{dashkovskiy2019stability}. Finally, we characterize the UGAS property of the induced system in terms of small-gain conditions and give useful sufficient conditions for it.%

\textbf{New vistas.} Our nonlinear ISS small-gain theorem (Theorem~\ref{thm_mainres}) has been developed for the case of the maximum formulation of the ISS property for the subsystems. However, our proof technique can be used for other formulations of the ISS property (semimaximum, summation, etc.), leading to very flexible small-gain results. In contrast to that, our method for the construction of the paths of strict decay for the gain operator relies strongly on the maximum formulation of the ISS property for subsystems. If the internal gains are all linear, and the ISS property for subsystems is given in a sum formulation, then the gain operator is linear and a path of strict decay with respect to such a gain operator exists (and can be explicitly constructed) if and only if the spectral radius of $\Gamma$ is less than one. Numerous characterizations of this condition have been developed in \cite{GlM21}. Some of these characterizations, and in particular, the method for the construction of paths of strict decay, have been extended in \cite{homogeneous} to the case of homogeneous and subadditive gain operators. The general nonlinear case remains, however, a challenging open problem.%

\textbf{Notation.} We write $\R$ for the reals and $\Z$ for the integers. The sets of nonnegative reals and integers are denoted by $\R_+$ and $\Z_+$, respectively. By $C^0(X,Y)$, we denote the set of all continuous mappings from a space $X$ to a space $Y$. In any metric space, we write $B_{\delta}(x)$ for the open ball of radius $\delta>0$ centered at $x$, and $\inner(A)$ for the interior of a subset $A \subset X$. We use the following classes of comparison functions:%
\begin{align*}
  \PC &:= \left\{\gamma \in C^0(\R_+,\R_+): \gamma(0) = 0,\ \gamma(r) > 0,\ \forall r > 0 \right\},\\
  \K &:= \left\{\gamma \in \PC: \ \gamma\mbox{ is strictly increasing}\right\}, \\
  \K_{\infty} &:= \left\{\gamma\in\K:\ \gamma\mbox{ is unbounded}\right\}, \\
  \LL &:= \big\{\gamma \in C^0(\R_+,\R_+):\ \gamma\mbox{ is strictly decreasing with}\\
	                &\qquad\qquad\qquad\qquad\qquad \lim_{t\rightarrow\infty}\gamma(t)=0 \big\}, \\
  \KL &:= \{\beta \in C^0(\R_+^2,\R_+): \beta(\cdot,t)\in\K,\ \forall t \geq 0,\ \\
		&\qquad\qquad\qquad\qquad\qquad\qquad \beta(r,\cdot)\in {\LL},\ \forall r > 0\}.%
\end{align*}

We write $\ell_{\infty}$ for the space of bounded real sequences $s = (s_i)_{i\in\N}$, which is a Banach space with the norm $\|s\|_{\ell_{\infty}} := \sup_{i\in\N}|s_i|$. The \emph{positive cone} in $\ell_{\infty}$ is given by $\ell_{\infty}^+ := \{ s \in \ell_{\infty} : s_i \geq 0,\ \forall i \in \N \}$. We define $\unit := (1,1,1,\ldots) \in \ell_{\infty}^+$. By $e_i$, $i\in \N$, we denote the $i$-th unit vector in $\ell_{\infty}$. Given $s^1,s^2 \in \ell_{\infty}$, we write $s^1 \oplus s^2$ for the vector given by the componentwise maximum of $s^1$ and $s^2$. If $\N$ is replaced by another index set $I$, we also write $\ell_{\infty}(I)$ for the corresponding Banach space and $\ell_{\infty}^+(I)$ for its positive cone. The notation $L_{\infty}(\R_+,U)$ is used for the Banach space of essentially bounded strongly measurable functions from $\R_+$ into a Banach space $U$.%

A function $\lambda:\R_+ \rightarrow X$ into some space $X$ is called \emph{piecewise right-continuous} if there is a partition of $\R_+$ into disjoint subintervals, $\R_+ = [0,t_1) \cup [t_1,t_2) \cup [t_2,t_3) \cup \ldots$, such that the restriction of $\lambda$ to each of the subintervals is continuous.%

\section{Technical setup}

\subsection{Interconnections}

Consider a family of control systems of the form%
\begin{equation}\label{eq_subsystem}
  \Sigma_i:\quad \dot{x}_i = f_i(x_i,\bar{x}_i,u_i),\quad i \in \N.%
\end{equation}
This family comes with sequences $(n_i)_{i\in\N}$ and $(m_i)_{i\in\N}$ of positive integers as well as \emph{finite} (possibly empty) sets $I_i \subset \N$, $i \notin I_i$, such that the following assumptions are satisfied:%
\begin{itemize}
\item The \emph{state vector} $x_i$ is an element of $\R^{n_i}$.%
\item The \emph{internal input vector} $\bar{x}_i$ is composed of the state vectors $x_j$, $j \in I_i$, and thus is an element of $\R^{N_i}$, where $N_i := \sum_{j \in I_i}n_j$.%
\item The \emph{external input vector} $u_i$ is an element of $\R^{m_i}$.%
\item The \emph{right-hand side} $f_i:\R^{n_i} \tm \R^{N_i} \tm \R^{m_i} \rightarrow \R^{n_i}$ is a continuous function.%
\item For every initial state $x_{i0} \in \R^{n_i}$ and all essentially bounded inputs $\bar{x}_i(\cdot)$ and $u_i(\cdot)$, there is a unique solution of $\Sigma_i$, which we denote by $\phi_i(t,x_{i0},\bar{x}_i,u_i)$ (it may be defined only on a bounded time interval).%
\end{itemize}
For each $i\in\N$, we fix norms on the spaces $\R^{n_i}$ and $\R^{m_i}$, respectively (these norms can be chosen arbitrarily). For brevity in notation, we avoid adding an index to these norms, indicating to which space they belong, and simply write $|\cdot|$ for each of them. The interconnection of the systems $\Sigma_i$, $i\in\N$, is defined on the state space $X := \ell_{\infty}(\N,(n_i))$, where%
\begin{align*}
  \ell_{\infty}(\N,(n_i)) := \{ x = (x_i)_{i\in\N} : x_i \in \R^{n_i},\ \sup_{i \in \N}|x_i| < \infty \}.%
\end{align*}
This space is a Banach space with the $\ell_{\infty}$-type norm%
\begin{equation*}
  \|x\|_X := \sup_{i\in\N}|x_i|.%
\end{equation*}
The space of admissible external input values is likewise defined as the Banach space%
\begin{equation*}
  U := \ell_{\infty}(\N,(m_i)),\quad \|u\|_U := \sup_{i\in\N}|u_i|.%
\end{equation*}
We choose the class of admissible external input functions as%
\begin{equation*}
  \UC := \{u \in L_{\infty}(\R_+,U) : u \mbox{ is piecewise right-continuous} \},%
\end{equation*}
which will be equipped with the $L_{\infty}$-norm%
\begin{equation*}
  \|u\|_{\UC} := \esssup_{t\in\R_+} |u(t)|_U.%
\end{equation*}
We define the right-hand side of the interconnected system by%
\begin{equation*}
  f:X \tm U \rightarrow \prod_{i\in\N}\R^{n_i},\quad  f(x,u) := (f_i(x_i,\bar{x}_i,u_i))_{i\in\N}.%
\end{equation*}
Hence, the interconnected system can formally be written as the differential equation%
\begin{equation*}
  \Sigma:\quad \dot{x} = f(x,u).%
\end{equation*}
To make sense of this equation, we need to define an appropriate notion of solution. For fixed $(u,x^0) \in \UC \tm X$, a function $\lambda:J \rightarrow X$, where $J \subset \R$ is an interval of the form $[0,T)$ with $0 < T \leq \infty$, is called a \emph{solution} of the Cauchy problem%
\begin{equation*}
  \dot{x} = f(x,u),\quad x(0) = x^0,%
\end{equation*}
provided that $s \mapsto f(\lambda(s),u(s))$ is a locally integrable $X$-valued function (in the Bochner integral sense) and%
\begin{equation*}
  \lambda(t) = x^0 + \int_0^t f(\lambda(s),u(s))\, \rmd s \mbox{\quad for all\ } t \in J.%
\end{equation*}

We say that the system $\Sigma$ is \emph{well-posed} if for every initial state $x^0 \in X$ and every external input $u \in \UC$ there exists a unique maximal solution. We denote this solution by $\phi(\cdot,x^0,u):[0,t_{\max}(x^0,u)) \rightarrow X$, where $0 < t_{\max}(x^0,u) \leq \infty$.%

The following theorem provides sufficient conditions for well-posedness of $\Sigma$, see \cite[Cor.~III.3]{KMS21}.%

\begin{theorem}\label{thm_sigma_wd}
Consider the coupled system $\Sigma$, composed of the subsystems $\Sigma_i$, and let the following assumptions hold:%
\begin{enumerate}
\item[(i)] $f(x,u) \in X$ for all $(x,u) \in X \tm U$.%
\item[(ii)] $f(\cdot,u):X \rightarrow X$ is continuous for each $u\in U$.%
\item[(iii)] $f(x,\cdot):U \rightarrow X$ is continuous for each $x \in X$.%
\item[(iv)] For each $u \in \UC$ and $x^0 \in X$, there are $\delta>0$ and locally integrable functions $\ell,\ell_0:\R_+ \rightarrow \R_+$ such that%
\begin{align*}
  \|f(x^1,u(t)) - f(x^2,u(t))\|_X &\leq \ell(t)\|x^1 - x^2\|_X, \\
	  \|f(x^0,u(t))\|_X &\leq \ell_0(t), %
\end{align*}
for all $x^1,x^2 \in B_{\delta}(x^0)$ and almost all $t \in \R_+$.
\end{enumerate}
Then $\Sigma$ is well-posed.%
\end{theorem}

If $\Sigma$ is well-posed, one has%
\begin{equation}\label{eq_relation_global_local_sol}
  \pi_i(\phi(t,x^0,u)) = \phi_i(t,x^0_i,\bar{x}_i,u_i)%
\end{equation}
for all $t \in [0,t_{\max}(x^0,u))$ and $i \in \N$, where $\pi_i:X \rightarrow \R^{n_i}$ denotes the canonical projection onto the $i$-th component, $\bar{x}_i(\cdot) = (\pi_j(\phi(\cdot,x^0,u)))_{j\in I_i}$, and $x_i^0,u_i$ denote the $i$-th components of $x^0$ and $u$, respectively, see \cite[Sec.~3]{KMS21}.%

In the rest of the paper, we assume that the following holds.%

%\mir{In the next assumption and in the paper in general: solutions of equations, or solutions of systems?}

%\ck{I think, we can use both simultaneously.}

\begin{assumption}\label{ass:Well-posedness-BIC}
The system $\Sigma$ is well-posed, and all of its uniformly bounded maximal solutions $\phi(\cdot,x,u)$ are global, i.e., they exist on $\R_+$ (this latter property is also called boundedness-implies-continuation (BIC) property).
\end{assumption}

\begin{remark}\label{rem:BIC-property}
If the function $f$ is uniformly bounded on bounded balls, and Lipschitz continuous on bounded balls with respect to the first argument, then $\Sigma$ is well-posed, and for any $R>0$ there is $\tau_R>0$ such that for all $x$ with $\|x\|_X \leq R$ and $u$ with $\|u\|_{\UC}\leq R$ the solution $\phi(\cdot,x,u)$ exists at least on $[0,\tau_R]$, which easily implies the BIC property; see \cite[Thm.~4.3.4]{CaH98} for a related result for systems without inputs.% 
\end{remark}

\subsection{Input-to-state stability}

We now recall the definition of input-to-state stability.%

\begin{definition}\label{def_ISS}
A well-posed system $\Sigma$ is called \emph{(uniformly) input-to-state stable (ISS)} if it is forward complete and there exist $\beta \in \KC\LC$ and $\gamma \in \KC_{\infty}$ such that%
\begin {equation*}
  \|\phi(t,x,u)\|_X \leq \beta(\|x\|_X,t) + \gamma(\|u\|_{\UC})%
\end{equation*}
for all $(t,x,u) \in \R_+ \tm X \tm \UC$.%
\end{definition}

Input-to-state stability is most often verified via the construction of an ISS Lyapunov function which is defined as follows.%

\begin{definition}
A function $V:X \rightarrow \R_+$ is called an \emph{ISS Lyapunov function (in an implication form)} for $\Sigma$ if it satisfies the following properties:%
\begin{enumerate}
\item[(i)] $V$ is continuous.%
\item[(ii)] There exist $\psi_1,\psi_2 \in \KC_{\infty}$ such that%
\begin{equation}\label{eq_isslf_coercivity}
  \psi_1(\|x\|_X) \leq V(x) \leq \psi_2(\|x\|_X) \mbox{\quad for all\ } x \in X.%
\end{equation}
\item[(iii)] There exist $\gamma\in\KC$ and $\alpha\in\PC$ such that for all $x\in X$ and $u\in\UC$ the following implication holds:%
\begin{equation}\label{eq_Lyap_impl}
  V(x) > \gamma(\|u\|_{\UC}) \quad \Rightarrow \quad \rmD^+ V_u(x) \leq -\alpha(V(x)),%
\end{equation}
where $\rmD^+ V_u(x)$ denotes the right upper Dini orbital derivative, defined as%
\begin{equation*}
  \rmD^+ V_u(x) := \limsup_{t \rightarrow 0^+} \frac{V(\phi(t,x,u)) - V(x)}{t}.%
\end{equation*}
\end{enumerate}
\end{definition}

The importance of ISS Lyapunov functions is due to the following result (cf.~\cite[Thm.~2.17]{MiP20}).%

\begin{proposition}\label{prop:Direct-Lyapunov-Theorem} 
If $\Sigma$ admits an ISS Lyapunov function, then it is ISS.
\end{proposition}

The construction of an ISS Lyapunov function is a complex problem, which becomes especially challenging if the system is nonlinear and of large size. In this paper, we assume that all components $\Sigma_i$ of an infinite network are ISS with corresponding ISS Lyapunov functions $V_i$. To find an ISS Lyapunov function $V$ for $\Sigma$, we exploit the interconnection structure and construct $V$ from $V_i$. Hence, we make the following assumption.%

\begin{assumption}\label{ass_subsystem_iss}
For each $i \in \N$, there exists a continuous function $V_i:\R^{n_i} \rightarrow \R_+$ which is continuously differentiable outside of $x_i = 0$ and satisfies the following properties:%
\begin{enumerate}
\item[(L1)] There exist $\psi_{i1},\psi_{i2} \in \KC_{\infty}$ such that%
\begin{equation}\label{eq_subsystem_iss_coerc}
  \psi_{i1}(|x_i|) \leq V_i(x_i) \leq \psi_{i2}(|x_i|) \mbox{\quad for all\ } x_i \in \R^{n_i}.%
\end{equation}
\item[(L2)] There exist $\gamma_{ij} \in \KC \cup \{0\}$, where $\gamma_{ij} = 0$ for all $j \in \N \setminus I_i$, and $\gamma_{iu} \in \KC$ as well as $\alpha_i \in \PC$ such that for all $x = (x_j)_{j\in\N} \in X$ and $u = (u_j)_{j\in\N} \in U$ the following implication holds:%
\begin{align}\label{eq_subsystem_orbitalder_est}
\begin{split}
  V_i(x_i) &> \max\Bigl\{\sup_{j\in I_i}\gamma_{ij}(V_j(x_j)),\gamma_{iu}(|u_i|) \Bigr\} \\
	&\Rightarrow \nabla V_i(x_i)f_i(x_i,\bar{x}_i,u_i) \leq -\alpha_i(V_i(x_i)).%
\end{split}
\end{align}
\end{enumerate}
The function $V_i$ is called an \emph{ISS Lyapunov function for $\Sigma_i$}. The functions $\gamma_{ij}$ and $\gamma_{iu}$ are called \emph{internal gains} and \emph{external gains}, respectively.%
\end{assumption}

Using the internal gains $\gamma_{ij}$ from Assumption \ref{ass_subsystem_iss}, we define the \emph{gain operator} $\Gamma:\ell_{\infty}^+ \rightarrow \ell_{\infty}^+$ by
\begin{equation}
\label{eq_gamma_def}
  \Gamma(s) := \Bigl( \sup_{j\in\N}\gamma_{ij}(s_j) \Bigr)_{i\in\N}.%
\end{equation}

In general, $\Gamma$ might be neither well-defined nor continuous. The following assumption guarantees both, see \cite[Lem.~2.1]{dashkovskiy2019stability} and \cite[Prop.~2]{MKG20}.%

\begin{assumption}\label{ass_gainop_wd}
The family $\{\gamma_{ij} : i,j \in \N\}$ is pointwise equicontinuous. That is, for every $r \geq 0$ and $\ep>0$, there exists $\delta = \delta(r,\ep) > 0$ such that $|r - \tilde{r}| \leq \delta$, $\tilde{r} \in \R_+$, implies $|\gamma_{ij}(r) - \gamma_{ij}(\tilde{r})| \leq \ep$ for all $i,j\in\N$.
\end{assumption}

Additionally, we make the following assumption on the external gains.%

\begin{assumption}\label{ass_gammaiu_max}
There is $\gamma_{\max}^u \in \KC$ such that $\gamma_{iu} \leq \gamma_{\max}^u$ for all $i \in \N$.%
\end{assumption}

We now introduce the concept of a \emph{path of strict decay} (for the gain operator $\Gamma$) which is of crucial importance in the construction of an ISS Lyapunov function for the interconnected system.%

\begin{definition}\label{def_omega_path}
A mapping $\sigma:\R_+ \rightarrow \ell_{\infty}^+$ is called a \emph{path of strict decay (for $\Gamma$)}, if the following properties hold:%
\begin{enumerate}
\item[(i)] There exists a function $\rho \in \KC_{\infty}$ such that%
\begin{equation*}
  \Gamma(\sigma(r)) \leq (\id + \rho)^{-1} \circ \sigma(r) \mbox{\quad for all\ } r \geq 0,%
\end{equation*}
where $(\id + \rho)^{-1}$ is applied componentwise.%
\item[(ii)] There exist $\sigma_{\min},\sigma_{\max} \in \KC_{\infty}$ satisfying%
\begin{equation*}
  \sigma_{\min} \leq \sigma_i \leq \sigma_{\max} \mbox{\quad for all\ } i \in \N.%
\end{equation*}
\item[(iii)] Each component function $\sigma_i = \pi_i \circ \sigma$, $i\in\N$, is a $\KC_{\infty}$-function.%
\item[(iv)] For every compact interval $K \subset (0,\infty)$, there exist $0 < c \leq C < \infty$ such that for all $r_1,r_2 \in K$ and $i \in \N$%
\begin{equation*}
  c|r_1 - r_2| \leq |\sigma_i^{-1}(r_1) - \sigma_i^{-1}(r_2)| \leq C|r_1 - r_2|.%
\end{equation*}
\end{enumerate}
\end{definition}

In Section~\ref{sec:Paths of strict decay}, we provide a method to construct paths of strict decay under suitable assumptions.%

\section{Nonlinear small-gain theorem}

Now we are able to present our small-gain result which yields an ISS Lyapunov function for the interconnected system $\Sigma$.%

\begin{theorem}
\label{thm_mainres}
Consider the interconnected system $\Sigma$, composed of subsystems $\Sigma_i$, $i\in\N$, and let the following assumptions hold.
\begin{enumerate}
\item[(i)] The system $\Sigma$ is well-posed and satisfies the BIC property (Assumption \ref{ass:Well-posedness-BIC}).
\item[(ii)] There exist ISS Lyapunov functions $V_i$ for the subsystems $\Sigma_i$ with associated internal gains $\gamma_{ij}$ and external gains $\gamma_{iu}$ (Assumption \ref{ass_subsystem_iss}). Moreover, there exist $\psi_1,\psi_2 \in \KC_{\infty}$ such that%
\begin{equation}\label{eq_coerc_unif_bounds}
  \psi_1 \leq \psi_{i1} \mbox{\quad and \quad} \psi_{i2} \leq \psi_2 \mbox{\quad for all\ } i \in \N.%
\end{equation}
\item[(iii)] The family $\{\gamma_{ij}\}$ of internal gains is pointwise equicontinuous and the external gains $\gamma_{iu}$ are uniformly upper bounded by a $\KC$-function (Assumptions \ref{ass_gainop_wd} and \ref{ass_gammaiu_max}).%
\item[(iv)] There exists a path $\sigma:\R_+ \rightarrow \ell_{\infty}^+$ of strict decay for the gain operator $\Gamma$, defined via the internal gains $\gamma_{ij}$.%
\item[(v)] For each $R>0$, there is a constant $L(R)>0$ such that%
\begin{equation}\label{eq_lyap_lipschitz}
  |V_i(x_i) - V_i(y_i)| \leq L(R)|x_i - y_i|%
\end{equation}
for all $i\in\N$ and $x_i,y_i \in B_R(0) \subset \R^{n_i}$.%
\item[(vi)] There exists $\tilde{\alpha} \in \PC$ such that $\alpha_i \geq \tilde{\alpha}$ for all $i\in\N$.%
\end{enumerate}
Then $\Sigma$ is ISS and an ISS Lyapunov function for $\Sigma$ is given by%
\begin{equation}\label{eq_def_V}
  V(x) := \sup_{i\in\N} \sigma_i^{-1}(V_i(x_i)) \mbox{\quad for all\ } x \in X.%
\end{equation}
Moreover, $V$ is locally Lipschitz continuous on $X \setminus \{0\}$. 
\end{theorem}

\begin{proof}
The proof proceeds in six steps.%

\emph{Step 1}: We show that $V$ satisfies inequalities of the form \eqref{eq_isslf_coercivity}, which also proves that $V$ assumes finite nonnegative values and $V(0)=0$. Using \eqref{eq_subsystem_iss_coerc} and \eqref{eq_coerc_unif_bounds}, we see that%
\begin{align*}
  \sigma_i^{-1}(V_i(x_i)) &\leq \sigma_{\min}^{-1}(\psi_{i2}(|x_i|)) \\
	&\leq \sigma_{\min}^{-1}(\psi_2(|x_i|)) \leq \sigma_{\min}^{-1}(\psi_2(\|x\|_X)).%
\end{align*}
Analogously, we obtain the lower estimate%
\begin{equation*}
  \sigma_i^{-1}(V_i(x_i)) \geq \sigma_{\max}^{-1}(\psi_{i1}(|x_i|)) \geq \sigma_{\max}^{-1}(\psi_1(|x_i|)).%
\end{equation*}
Together, these estimates imply%
\begin{align*}
  \sigma_{\max}^{-1} \circ \psi_1(\|x\|_X) &= \sup_{i\in\N} \sigma_{\max}^{-1}(\psi_1(|x_i|)) \\
	&\leq V(x) \leq \sigma_{\min}^{-1} \circ \psi_2(\|x\|_X).%
\end{align*}
Since $\sigma_{\max}^{-1} \circ \psi_1$ and $\sigma_{\min}^{-1} \circ \psi_2$ are $\KC_{\infty}$-functions, the desired coercivity estimates hold.%

\emph{Step 2}: We prove that $V$ is continuous and locally Lipschitz continuous outside of $x=0$. Continuity at $x=0$ follows from coercivity as shown in Step 1. Hence, let $0 \neq x \in X$. Define%
\begin{equation*}
  \delta = \delta(x) := \frac{1}{3}\psi_2^{-1} \circ \sigma_{\min} \circ \sigma_{\max}^{-1} \circ \psi_1\Bigl(\frac{\|x\|_X}{4}\Bigr) \leq \frac{\|x\|_X}{12}.%
\end{equation*}
For all $y \in B_{\delta}(x)$, it holds that $\|y\|_X \geq \frac{1}{2}\|x\|_X$, implying%
\begin{equation}\label{eq_view}
  V(y) \geq \sigma_{\max}^{-1} \circ \psi_1(\|y\|_X) \geq \sigma_{\max}^{-1} \circ \psi_1\Bigl(\frac{\|x\|_X}{2}\Bigr).%
\end{equation}
Define $I_{\delta} := \{i \in \N : |x_i| \geq 2\delta\}$. For any $y \in B_{\delta}(x)$ and $i \notin I_{\delta}$, we have%
\begin{equation*}
  |y_i| \leq |x_i| + |x_i - y_i| \leq 3\delta = \psi_2^{-1} \circ \sigma_{\min} \circ \sigma_{\max}^{-1} \circ \psi_1\Bigl(\frac{\|x\|_X}{4}\Bigr).%
\end{equation*}
This implies%
\begin{equation*}
  \sigma_i^{-1}(V_i(y_i)) \leq \sigma_{\min}^{-1} \circ \psi_2(|y_i|) \leq \sigma_{\max}^{-1} \circ \psi_1\Bigl(\frac{\|x\|_X}{4}\Bigr).%
\end{equation*}
In view of \eqref{eq_view}, we see that%
\begin{equation*}
  V(y) = \sup_{i \in I_{\delta}}\sigma_i^{-1}(V_i(y_i))%
\end{equation*}
for all $y \in B_{\delta}(x)$. It is then easy to see that%
\begin{equation*}
  |V(y^1) - V(y^2)| \leq \sup_{i \in I_{\delta}}|\sigma_i^{-1}(V_i(y^1_i)) - \sigma_i^{-1}(V_i(y^2_i))|%
\end{equation*}
for all $y^1,y^2 \in B_{\delta}(x)$. Since $V_i(y^1_i)$ and $V_i(y^2_i)$ on the right-hand side of this inequality are contained in the compact interval $[\psi_1(\delta),\psi_2(\|x\|_X + \delta)] \subset (0,\infty)$, the definition of a path of strict decay together with Assumption (v) implies the existence of a constant $L > 0$ such that for all $y^1,y^2\in B_{\delta}(x)$%
\begin{equation*}
  |V(y^1) - V(y^2)| \leq \sup_{i \in I_{\delta}}L |y^1_i - y^2_i| \leq L \|y^1 - y^2\|_X.%
\end{equation*}
This proves Lipschitz continuity of $V$ on $B_{\delta}(x)$.%

\emph{Step 3}: We prove an auxiliary result needed for the proof of implication \eqref{eq_Lyap_impl}. Let $\rho \in \KC_{\infty}$ satisfy%
\begin{equation}\label{eq_omegapath_cond}
  \Gamma(\sigma(r)) \leq (\id + \rho)^{-1} \circ \sigma(r) \mbox{\quad for all\ } r\in\R_+,%
\end{equation}
as required in the definition of a path of strict decay. Then we fix a function $\mu \in \KC_{\infty}$ such that $\mu(r) < \rho(r)$ for all $r > 0$ and introduce for every $0 \neq x \in X$ the set%
\begin{equation}\label{eq_i_choice}
  I(x) := \bigl\{ i \in \N : V(x) \leq \sigma_i^{-1}((\id + \mu)(V_i(x_i))) \bigr\}.%
\end{equation}
Now we prove the following claim:%

\begin{framed}
Every $x \neq 0$ has a neighborhood in $X$ on which all of the functions $\sigma_i^{-1} \circ V_i \circ \pi_i$, $i \in \N \setminus I(x)$, are bounded away from $V$, formally:%
\begin{align}\label{eq_ix_rep}
\begin{split}
  \exists \ep,\delta>0:\ & \|x - y\|_X \leq \delta \wedge i \in \N \setminus I(x) \\
	&\Rightarrow \ V(y) > \sigma_i^{-1}(V_i(y_i)) + \ep.%
\end{split}
\end{align}
In particular, this shows that $I(x) \neq \emptyset$.%
\end{framed}

Assume towards a contradiction that the claim is false. Then we can find sequences $y^n \rightarrow x$ and $i_n \in \N \setminus I(x)$ such that%
\begin{equation}\label{eq_contradictpos}
  V(y^n) \leq \sigma_{i_n}^{-1}(V_{i_n}(y^n_{i_n})) + \frac{1}{n} \mbox{\quad for all\ } n \in \N.%
\end{equation}
At the same time, $i_n \notin I(x)$ implies%
\begin{equation*}
  V(x) > \sigma_{i_n}^{-1}((\id + \mu)(V_{i_n}(x_{i_n}))) \mbox{\quad for all\ } n \in \N.%
\end{equation*}
Combining these two inequalities, we obtain%
\begin{equation*}
  V(x) - V(y^n) > \sigma_{i_n}^{-1}((\id + \mu)(V_{i_n}(x_{i_n}))) - \sigma_{i_n}^{-1}(V_{i_n}(y^n_{i_n})) - \frac{1}{n}.%
\end{equation*}
We can find a compact interval $K \subset (0,\infty)$ such that for sufficiently large $n$ we have $(\id + \mu)(V_{i_n}(x_{i_n})), V_{i_n}(y^n_{i_n}) \in K$. Indeed, this follows from the estimates%
\begin{itemize}
\item $(\id + \mu)(V_{i_n}(x_{i_n})) \leq (\id + \mu) \circ \psi_2(\|x\|_X)$;%
\item $V_{i_n}(y^n_{i_n}) \leq \psi_2(\|y^n\|_X) \leq \psi_2(2\|x\|_X)$ for all $n$ large enough;%
\item $V_{i_n}(y^n_{i_n}) \geq \sigma_{i_n}(V(y^n) - 1/n) \geq \sigma_{\min}(V(y^n) - 1/n) \geq \sigma_{\min}(\sigma_{\max}^{-1}(\psi_1(\frac{1}{2}\|x\|_X)) - 1/n) > 0$ for all $n$ large enough;%
\item $(\id + \mu)(V_{i_n}(x_{i_n})) \geq (\id + \mu) \circ \psi_1(|x_{i_n}|) \geq (\id + \mu) \circ \psi_1(|y^n_{i_n}| - |x_{i_n} - y^n_{i_n}|) \geq (\id + \mu) \circ \psi_1(\frac{1}{2}|y^n_{i_n}|) \geq (\id + \mu) \circ \psi_1(\frac{1}{2}\psi_2^{-1}(V_{i_n}(y^n_{i_n})))$ for all $n$ large enough; this can be lower bounded by using the previous estimates.
\end{itemize}
By the definition of a path of strict decay, we thus have constants $0 < c \leq C < \infty$ such that%
\begin{align*}
  &|\sigma_{i_n}^{-1}((\id+\mu)(V_{i_n}(x_{i_n}))) - \sigma_{i_n}^{-1}(V_{i_n}(y^n_{i_n}))| \\
	&\qquad = c_n|(\id+\mu)(V_{i_n}(x_{i_n})) - V_{i_n}(y^n_{i_n})|%
\end{align*}
for some $c_n \in [c,C]$, when $n$ is large enough, which (by monotonicity of $\sigma_{i_n}^{-1}$) yields%
\begin{equation*}
  V(x) - V(y^n) > c_n((\id + \mu)(V_{i_n}(x_{i_n})) - V_{i_n}(y^n_{i_n})) - \frac{1}{n}.%
\end{equation*}
From \eqref{eq_lyap_lipschitz} it follows that for some $L>0$ (independent of $n$),%
\begin{align*}
  V_{i_n}(y^n_{i_n}) &= V_{i_n}(x_{i_n}) + V_{i_n}(y^n_{i_n}) - V_{i_n}(x_{i_n}) \\
	&\leq V_{i_n}(x_{i_n}) + L \|y^n - x\|_X.%
\end{align*}
Putting $\delta_n := \|y^n - x\|_X$, we thus obtain%
\begin{align*}
  &V(x) - V(y^n) \\
	&\quad > c_n((\id + \mu)(V_{i_n}(x_{i_n})) - [V_{i_n}(x_{i_n}) + L\delta_n]) - \frac{1}{n} \\
	         &\quad = c_n\mu(V_{i_n}(x_{i_n})) - c_nL\delta_n - \frac{1}{n} \\
					&\quad \geq c_n \mu(\psi_1(|x_{i_n}|)) - c_nL\delta_n - \frac{1}{n}.%
\end{align*}
We can also write this as%
\begin{align*}
  & 0 \leq \mu \circ \psi_1(|x_{i_n}|) < b_n,\\
	 b_n &:= c_n^{-1}\Bigl([V(x) - V(y^n)] + c_nL\delta_n + \frac{1}{n}\Bigr).%
\end{align*}
Note that $b_n \rightarrow 0$ as $n \rightarrow \infty$, which implies $|x_{i_n}| \rightarrow 0$ as $n \rightarrow \infty$, and since $y^n \rightarrow x$, also $|y^n_{i_n}| \rightarrow 0$. Then \eqref{eq_contradictpos} yields%
\begin{equation*}
  0 \leq V(y^n) \leq \sigma_{\min}^{-1} \circ \psi_2(|y^n_{i_n}|) + \frac{1}{n} \rightarrow 0.%
\end{equation*}
Since $V(y^n) \rightarrow V(x)$, we obtain $V(x) = 0$, and hence $x = 0$, a contradiction.%

\emph{Step 4}: We define the $\KC$-function $\gamma$ by%
\begin{equation*}
  \gamma(r) := \sigma_{\min}^{-1} \circ (\id + \rho) \circ \gamma_{\max}^u(r) \mbox{\quad for all\ } r \geq 0,%
\end{equation*}
and prove the following claim:%

\begin{framed}
The inequality $V(x) > \gamma(\|u\|_{\UC})$ for some $x \in X$ and $u \in \UC$ implies the existence of $T>0$ such that%
\begin{equation}\label{eq_iss_subsystem_est}
  \nabla V_i(\phi_i(t)) f_i(\phi_i(t),\bar{\phi}_i(t),u_i(t)) \leq -\tilde{\alpha}(V_i(\phi_i(t)))%
\end{equation}
for all $i \in I(x)$ and $t \in [0,T]$, where $\phi(t) := \phi(t,x,u)$, $\phi_j(t)$ is the $j$-th component of $\phi(t)$ (for every $j\in\N$) and $\bar{\phi}_i(t) = (\phi_j(t))_{j\in I_i}$.%
\end{framed}

Let us fix $x$ and $u$ as in the claim and note that $x \neq 0$. As $V$ and $\phi(\cdot)$ are continuous, $V(x) > \gamma(\|u\|_{\UC})$ implies%
\begin{align}\label{eq_lyap_hypothesis}
\begin{split}
  V(\phi(t)) &> \sigma_{\min}^{-1} \circ (\id + \rho) \circ \gamma_{\max}^u(\|u\|_{\UC})\\
	&\geq \sigma_{\min}^{-1} \circ (\id + \rho) \circ \gamma_{iu}(|u_i(t)|)%
\end{split}
\end{align}
for all $i \in \N$ and $t \in [0,T]$, where $T>0$ is chosen sufficiently small. If $T$ is chosen further small enough, then%
\begin{equation}\label{eq_subclaim}
  V(\phi(t)) < \sigma_i^{-1} \circ (\id + \rho) \circ V_i(\phi_i(t))\ \forall i \in I(x),\ t \in [0,T].%
\end{equation}
We prove this by contradiction. Assume that there are sequences $0 < t_n \rightarrow 0$ and $i_n \in I(x)$ such that%
\begin{equation*}
  V(\phi(t_n)) \geq \sigma_{i_n}^{-1} \circ (\id + \rho) \circ V_{i_n}(\phi_{i_n}(t_n)).%
\end{equation*}
The left-hand side of this inequality converges to $V(x)$ as $n \rightarrow \infty$. Hence, we find a sequence $0 < \ep_n \rightarrow 0$ such that%
\begin{align*}
  & \sigma_{i_n}^{-1} \circ (\id + \rho) \circ V_{i_n}(\phi_{i_n}(t_n)) \leq V(x) + \ep_n \\
	&\qquad \leq \sigma_{i_n}^{-1} \circ (\id + \mu) \circ V_{i_n}(x_{i_n}) + \ep_n.%
\end{align*}
We can further estimate%
\begin{align*}
  V_{i_n}(\phi_{i_n}(t_n)) &\geq V_{i_n}(x_{i_n}) - |V_{i_n}(x_{i_n}) - V_{i_n}(\phi_{i_n}(t_n))| \\
	                         &\geq V_{i_n}(x_{i_n}) - L |x_{i_n} - \phi_{i_n}(t_n)| \\
													 &\geq V_{i_n}(x_{i_n}) - L\|x - \phi(t_n)\|_X,%
\end{align*}
where we use that $|\phi_{i_n}(t_n)| \leq |\phi_{i_n}(t_n) - x_{i_n}| + |x_{i_n}| \leq \|\phi(t_n) - x\|_X + \|x\|_X$ (implying the existence of $L>0$ by Assumption (v)). Hence, we obtain%
\begin{align*}
  &\sigma_{i_n}^{-1} \circ (\id + \rho)(V_{i_n}(x_{i_n}) - L\|x - \phi(t_n)\|_X) \\
	&\qquad \leq \sigma_{i_n}^{-1} \circ (\id + \mu)(V_{i_n}(x_{i_n})) + \ep_n.%
\end{align*}
We write this inequality as%
\begin{align*}
  \ep_n &\geq \sigma_{i_n}^{-1} \circ (\id + \rho)(V_{i_n}(x_{i_n}) - L\|x - \phi(t_n)\|_X) \\
	&\qquad - \sigma_{i_n}^{-1} \circ (\id + \mu)(V_{i_n}(x_{i_n})).%
\end{align*}
With a similar reasoning as used before, we can show that the arguments of $\sigma_{i_n}^{-1}$ are contained in a compact subset of $(0,\infty)$ for all sufficiently large $n$. Hence, there are numbers $c_n \in [c,C]$ such that%
\begin{align*}
  \ep_n &\geq c_n \Bigl[\rho(V_{i_n}(x_{i_n}) - L\|x - \phi(t_n)\|_X) \\
	&\qquad- \mu(V_{i_n}(x_{i_n})) - L\|x - \phi(t_n)\|_X\Bigr].%
\end{align*}
Using that $V_{i_n}(x_{i_n})$ is contained in a compact interval for all $n$ and $\rho$ is uniformly continuous on this interval, we find a sequence $0 < \delta_n \rightarrow 0$ such that%
\begin{equation*}
  \rho(V_{i_n}(x_{i_n}) - L\|x - \phi(t_n)\|_X) \geq \rho(V_{i_n}(x_{i_n})) - \delta_n\ \forall n \in \N.%
\end{equation*}
Also using that $\rho - \mu > 0$ (on $(0,\infty)$) and $V_{i_n}(x_{i_n}) \geq (\id + \mu)^{-1} \circ \sigma_{\min}(V(x))$, this implies%
\begin{align*}
  \ep_n &\geq c_n \Bigl[ (\rho - \mu)( (\id + \mu)^{-1} \circ \sigma_{\min}(V(x))) \\
	&\qquad\qquad - \delta_n - L\|x - \phi(t_n)\|_X \Bigr].%
\end{align*}
Hence, as $n \rightarrow \infty$, we obtain the contradiction $V(x) = 0$, as $V(x) > \gamma(\|u\|_{\UC})\geq 0$. This proves \eqref{eq_subclaim}.%

From \eqref{eq_subclaim}, it then follows that for all $i \in I(x)$ and $t \in [0,T]$:%
\begin{align*}
  V_i(\phi_i(t)) &> (\id + \rho)^{-1} \circ \sigma_i(V(\phi(t))) \stackrel{\eqref{eq_omegapath_cond}}{\geq} \Gamma_i(\sigma(V(\phi(t)))) \\
	&= \sup_{j\in I_i} \gamma_{ij}(\sigma_j(V(\phi(t)))) \geq \sup_{j \in I_i} \gamma_{ij}(V_j(\phi_j(t))).%
\end{align*}
At the same time, \eqref{eq_lyap_hypothesis} together with \eqref{eq_subclaim} implies%
\begin{align*}
  V_i(\phi_i(t)) &> (\id + \rho)^{-1} \circ \sigma_i \circ \sigma_{\min}^{-1} \circ (\id + \rho) \circ \gamma_{iu}(|u_i(t)|) \\
	&\geq \gamma_{iu}(|u_i(t)|).%
\end{align*}
Putting both estimates together, we obtain%
\begin{equation*}
  V_i(\phi_i(t)) > \max\Bigl\{ \sup_{j \in I_i}\gamma_{ij}(V_j(\phi_j(t))),\gamma_{iu}(|u_i(t)|)\Bigr\}.%
\end{equation*}
By \eqref{eq_subsystem_orbitalder_est} together with Assumption (vi), this implies \eqref{eq_iss_subsystem_est}, which proves the claim.%

\emph{Step 5}: We show the implication \eqref{eq_Lyap_impl}. First observe that the case $x=0$ does not occur, since $V(0) = 0 > \gamma(\|u\|_{\UC})$ is never satisfied. Hence, let us fix $0 \neq x \in X$ and $u \in \UC$ satisfying $V(x) > \gamma(\|u\|_{\UC})$. By Step 3, there is $\delta>0$ such that%
\begin{equation*}
  V(y) = \sup_{i \in I(x)} \sigma_i^{-1}(V_i(y_i)) \mbox{\quad for all\ } y \in B_{\delta}(x).%
\end{equation*}
Step 4 shows that for all sufficiently small $t\geq0$ we have%
\begin{equation*}
  \nabla V_i(\phi_i(t)) f_i(\phi_i(t),\bar{\phi}_i(t),u_i(t)) \leq -\tilde{\alpha}(V_i(\phi_i(t)))\ \forall i \in I(x).%
\end{equation*}
Now let us introduce the Cauchy problem%
\begin{equation*}
  \dot{v}(t) = -\tilde{\alpha}(v(t)),\quad v(0) = v_0 \in \R_+.%
\end{equation*}
By \cite[Lem.~6]{MKG20}, we can assume that $\tilde{\alpha}$ is globally Lipschitz. Then the Cauchy problem has a globally defined unique solution that we denote by $\VC(t,v_0)$. For every $i \in I(x)$ and all sufficiently small $t\geq0$, Lemma \ref{lem_techlem1} guarantees that%
\begin{equation*}
  V_i(\phi_i(t)) \leq \VC(t,V_i(x_i)).%
\end{equation*}
It thus follows that%
\begin{align*}
 \frac{1}{t}(V(\phi(t)) &- V(x)) \\
	&= \frac{1}{t}\Bigl[\sup_{i\in I(x)} \sigma_i^{-1}(V_i(\phi_i(t))) - \sup_{i \in I(x)} \sigma_i^{-1}(V_i(x_i))\Bigr] \\
	&\leq\frac{1}{t}\sup_{i\in I(x)}\bigl[\sigma_i^{-1}(V_i(\phi_i(t))) - \sigma_i^{-1}(V_i(x_i))\bigr] \\
	&\leq \frac{1}{t}\sup_{i\in I(x)}\bigl[\sigma_i^{-1}(\VC(t,V_i(x_i))) - \sigma_i^{-1}(V_i(x_i))\bigr] \\
	&= \frac{1}{t}\sup_{i\in I(x)} -\bigl|\sigma_i^{-1}(\VC(t,V_i(x_i))) - \sigma_i^{-1}(V_i(x_i))\bigr|.%
\end{align*}
For all $i \in I(x)$ and all $t\geq 0$, we have%
\begin{equation*}
  \VC(t,V_i(x_i)) \leq V_i(x_i) \leq \sigma_i \circ V(x) \leq \sigma_{\max}(V(x))%
\end{equation*}
and%
\begin{equation*}
  V_i(x_i) \geq (\id+\mu)^{-1}\circ\sigma_i \circ V(x) \geq (\id+\mu)^{-1}\circ\sigma_{\min} \circ V(x).
\end{equation*}
With $t^*>0$ chosen small enough, for all $t \in(0,t^*)$, we have%  
\begin{align*}
  \VC(t,V_i(x_i)) &\geq \VC(t,(\id+\mu)^{-1}\circ\sigma_{\min} \circ V(x)) \\
	&\geq \frac{1}{2} (\id+\mu)^{-1}\circ\sigma_{\min} \circ V(x).%
\end{align*}
Now, define%
\begin{equation*}
  K(r) := \Big[\frac{1}{2}(\id+\mu)^{-1}(\sigma_{\min}(r)),\sigma_{\max}(r)\Big]%
\end{equation*}
and let $c = c(K(r)) > 0$ be the maximal constant such that%
\begin{equation*}
  |\sigma_i^{-1}(r_1) - \sigma_i^{-1}(r_2)| \geq c|r_1 - r_2|, \quad \forall r_1,r_2 \in K(r).%
\end{equation*}
For all $t\in(0,t^*)$, we obtain% 
\begin{align*}
  & \frac{1}{t}(V(\phi(t)) - V(x)) \\
	&\leq -c(K(V(x))) \inf_{i\in I(x)} \frac{1}{t} \bigl(V_i(x_i) - \VC(t,V_i(x_i))\bigr) \\
	&= -c(K(V(x))) \inf_{i\in I(x)} \frac{1}{t} \int_0^t \tilde{\alpha}(\VC(s,V_i(x_i)))\, \rmd s \\
	&\leq -c(K(V(x))) \frac{1}{t}\int_0^t \min_{\varrho \in K(V(x))} \tilde{\alpha}(\VC(s,\varrho))\, \rmd s.%
\end{align*}
Observe that the function%
\begin{equation*}
  s \mapsto \min_{\varrho \in K(V(x))} \tilde{\alpha}(\VC(s,\varrho))%
\end{equation*}
is continuous as $\VC(\cdot,\cdot)$ is continuous, and thus uniformly continuous on compact sets. Hence,%
\begin{equation*}
  \rmD^+ V_u(x) \leq -c(K(V(x))) \min_{\varrho \in K(V(x))}\tilde{\alpha}(\varrho).%
\end{equation*}
Therefore, we have proved the implication%
\begin{equation*}
  V(x) > \gamma(\|u\|_{\UC}) \quad \Rightarrow \quad \rmD^+ V_u(x) \leq -\hat{\alpha}(V(x))%
\end{equation*}
for all $0 \neq x \in X$ with%
\begin{equation*}
  \hat{\alpha}(r) := c(K(r)) \min_{\varrho \in K(r)} \tilde{\alpha}(\varrho),\quad \forall r > 0.%
\end{equation*}

\emph{Step 6}: It remains to lower bound $\hat{\alpha}$ by a positive definite function. For each $r>0$, define $K_2(r) := \bigcup_{q \in K(r)}\mathrm{str}(1,q)$, where $\mathrm{str}(1,q)$ equals $[1,q]$ for $q\geq 1$ and $[q,1]$ for $q<1$. Clearly, $K_2(r)$ is a compact subset of $(0,\infty)$. Further, we introduce%
\begin{equation*}
  \hat{\alpha}_2(r) := c(K_2(r)) \min_{\varrho \in K_2(r)} \tilde{\alpha}(\varrho),\quad \forall r > 0.%
\end{equation*}
As $K(r)\subset K_2(r)$ for any $r>0$, $\hat{\alpha}(r) \geq \hat{\alpha}_2(r)$ for all $r>0$. Furthermore, there is $r_{\min}$ such that $K_2(r_1) \supset K_2(r_2)$ for all $r_1,r_2 \in (0,r_{\min})$ with $r_1 < r_2$. This implies that $\hat{\alpha}_2$ is a non-decreasing positive function on $(0,r_{\min})$, and $\lim_{r\downarrow0}\hat{\alpha}_2(r) = 0$. Moreover, for all $r\in (0,r_{\min})$ we have%
\begin{equation*}
  \hat{\alpha}_2(r) = \frac{2}{r}\int_{r/2}^r \hat{\alpha}_2(r)\, \rmd s \geq \frac{2}{r}\int_{r/2}^r \hat{\alpha}_2(s)\, \rmd s,%
\end{equation*}
where $\hat{\alpha}_2$ is integrable on $(0,r_{\min})$ as it is monotone on this interval. Hence, $\hat{\alpha}_2$ and thus $\hat{\alpha}$ can be lower bounded by a continuous function on $[0,r_{\min}]$. Similarly, there is $r_{\max}$ such that $K_2(r_1)\supset K_2(r_2)$ for all $r_1,r_2 \in (r_{\max},\infty)$ with $r_1 > r_2$. This implies that $\hat{\alpha}_2$ is a non-increasing positive function on $(r_{\max},\infty)$. Consequently, for $r\in (r_{\max},\infty)$, we have%
\begin{equation*}
  \hat{\alpha}_2(r) = \frac{1}{r}\int_{r}^{2r} \hat{\alpha}_2(r)\, \rmd s \geq \frac{1}{r}\int_{r}^{2r} \hat{\alpha}_2(s)\, \rmd s.%
\end{equation*}
Hence, $\hat{\alpha}_2$, and thus $\hat{\alpha}$, can be lower bounded by a continuous function on $[r_{\max},\infty)$. As $\hat{\alpha}$ assumes positive values and is bounded away from zero on every compact interval in $(0,\infty)$, $\hat{\alpha}$ can be lower bounded by a positive definite function, which we denote by $\alpha$, on $\R_+$. This implies%
\begin{equation}\label{eq_limp}
  V(x) > \gamma(\|u\|_{\UC}) \ \Rightarrow \ \rmD^+ V_u(x) \leq -\alpha(V(x))\ \forall x \in X.%
\end{equation}
Thus, $V$ is an ISS Lyapunov function for $\Sigma$. Proposition~\ref{prop:Direct-Lyapunov-Theorem} implies that $\Sigma$ is ISS. The proof is complete.%
\end{proof}

\section{Gain operators and their properties}\label{sec:Gain operators and their properties}

A crucial assumption in our small-gain theorem is the existence of a path of strict decay for the operator $\Gamma$. Our next goal is to understand under which conditions such a path exists, and to provide an explicit expression for it. We base our analysis on the properties of the gain operator, derived in this section.%

From now on, we always assume that the family $\{\gamma_{ij}\}$ is pointwise equicontinuous (Assumption~\ref{ass_gainop_wd}), implying that $\Gamma$ is well-defined and continuous. We assume throughout this section that for each $i \in \N$, $\gamma_{ij} \neq 0$ only for finitely many $j\in\N$, though many of the following results hold without this assumption.%

The most important property of $\Gamma$ is its monotonicity: for all $s^1,s^2 \in \ell_{\infty}^+$ we have the implication%
\begin{equation*}
  s^1 \leq s^2 \quad \Rightarrow \quad \Gamma(s^1) \leq \Gamma(s^2).%
\end{equation*}
Moreover, we note that $\Gamma(0) = 0$ and that $\Gamma$ is a \emph{max-preserving operator}, i.e.%
\begin{equation*}
  \Gamma(s^1 \oplus s^2) = \Gamma(s^1) \oplus \Gamma(s^2) \mbox{\quad for all\ } s^1,s^2 \in \ell_{\infty}^+.%
\end{equation*}

Now, we recall the important \emph{robust and robust strong small-gain conditions}, introduced in \cite{dashkovskiy2019stability}, which are closely related to the stability properties of the discrete-time system induced by $\Gamma$. We modify these properties to make them more compatible with max-type gain operators.%

\begin{definition}\label{def:Small-gain-conditions}
We say that the operator $\Gamma$ satisfies%
\begin{enumerate}
\item[(i)] the \emph{small-gain condition (SGC)} if%
\begin{equation}\label{eq_sg_cond}
  \Gamma(s) \ngeq s \mbox{\quad for all\ } s \in \ell_{\infty}^+ \setminus \{0\}.%
\end{equation}
\item[(ii)] the \emph{strong small-gain condition} if there is $\rho \in \KC_{\infty}$ with%
\begin{equation}\label{eq_ssg_cond}
  D_{\rho} \circ \Gamma(s) \ngeq s \mbox{\quad for all\ } s \in \ell_{\infty}^+ \setminus \{0\}%
\end{equation}
for the operator $D_{\rho}:\ell_{\infty}^+ \rightarrow \ell_{\infty}^+$, defined by%
\begin{equation*}
  D_{\rho}(s) := \bigl((\id + \rho)(s_i)\bigr)_{i\in\N}.%
\end{equation*}
\item[(iii)] the \emph{max-robust small-gain condition} if there is $\omega \in \KC_{\infty}$ with $\omega<\id$ such that for all $i,j \in \N$ the operator%
\begin{equation}\label{eq_def_gammaij}
  \Gamma_{ij}(s) := \Gamma(s) \oplus \omega(s_j)e_i \mbox{\quad for all\ } s \in \ell_{\infty}^+%
\end{equation}
satisfies the small-gain condition.
\item[(iv)] the \emph{max-robust strong small-gain condition} if there are $\omega \in \KC_{\infty}$ with $\omega<\id$ and $\rho \in \KC_{\infty}$ such that for all $i,j\in\N$ the operator $\Gamma_{ij}$, defined in \eqref{eq_def_gammaij}, satisfies the strong small-gain condition with the same $\rho$ for all $i,j$.
\end{enumerate}
\end{definition}

In the next lemma, we introduce the so-called \emph{strong transitive closure} (or Kleene star operator) $Q$ of the gain operator $\Gamma$, which provides the crucial tool for the construction of a path of strict decay. This result was first shown in \cite[Lem.~4.3]{dashkovskiy2019stability}, strengthened in \cite[Lem.~B.5]{MKG20}, and is now even more strengthened, since the robust SGC is replaced by the less restrictive max-robust SGC.%

\begin{lemma}\label{lem_q_op}%
Assume that $\Gamma$ satisfies the max-robust SGC with some $\omega \in \KC_{\infty}$. Then the operator%
\begin{equation}\label{eq:Q-operator}
  Q(s) := \bigoplus_{k\in\Z_+}\Gamma^k(s) \mbox{\quad for all\ } s \in \ell_{\infty}^+,%
\end{equation}
is well-defined and has the following properties:%
\begin{equation}\label{eq:upper estimate on Q}
  s \leq Q(s) \leq \omega^{-1}(\|s\|_{\ell_{\infty}})\unit \mbox{\quad for all\ } s \in \ell_{\infty}^+,%
\end{equation}
\begin{equation}\label{eq_q_decay}
  \Gamma(Q(s)) \leq Q(s) \mbox{\quad for all\ } s \in \ell_{\infty}^+.%
\end{equation}
\end{lemma}

\begin{proof}
The proof is only a slight variation of the proof of \cite[Lem.~B.5]{MKG20}. Here, it suffices to show that the assumption $\sup_{k\in\Z_+}\Gamma^k_i(s) > \omega^{-1}(\|s\|_{\ell_{\infty}})$ for some $i \in \N$ and $s \in \ell_{\infty}^+$ leads to a contradiction. From this assumption, the existence of $j \in \N$ and $j_1,\ldots,j_{k-1} \in \N$ follows such that%
\begin{equation*}
  \gamma_{ij_1} \circ \gamma_{j_1j_2} \circ \cdots \circ \gamma_{j_{k-1}j}(s_j) > \omega^{-1}(\|s\|_{\ell_{\infty}}).%
\end{equation*}
Now, consider the operator $\Gamma_{ji}$, as defined in \eqref{eq_def_gammaij}:%
\begin{align*}
  \Gamma_{ji}(s) &= \Bigl( \max\{ \sup_{k\in\N} \gamma_{lk}(s_k), \omega(s_i) \delta_{jl} \} \Bigr)_{l \in \N} \\
	               &= \Bigl( \max\{ \sup_{k\in\N} \gamma_{lk}(s_k), \sup_{k\in\N} \omega(s_k) \delta_{jl}\delta_{ik} \} \Bigr)_{l \in \N} \\
								 &= \Bigl( \sup_{k\in\N} \max\{ \gamma_{lk}(s_k), \omega(s_k) \delta_{jl}\delta_{ik} \} \Bigr)_{l \in \N}.%
\end{align*}
Hence, $\Gamma_{ji}$ is a gain operator induced by the gains%
\begin{equation*}
  \tilde{\gamma}_{lk}(r) := \max\{\gamma_{lk}(r),\delta_{jl}\delta_{ik}\omega(r)\}.%
\end{equation*}
Since $\Gamma_{ji}$ satisfies the SGC by assumption, by \cite[Lem.~B.3]{MKG20}, all cycles built from the gains $\tilde{\gamma}_{lk}$ are contractions. This implies a contradiction, namely%
\begin{align*}
  s_j &> \tilde{\gamma}_{ji} \circ \tilde{\gamma}_{ij_1} \circ \tilde{\gamma}_{j_1j_2} \circ \cdots \circ \tilde{\gamma}_{j_{k-1}j}(s_j) \\
	&\geq \tilde{\gamma}_{ji} \circ \gamma_{ij_1} \circ \gamma_{j_1j_2} \circ \cdots \circ \gamma_{j_{k-1}j}(s_j) \\
	&= \max\{\gamma_{ji},\omega\} \circ \gamma_{ij_1} \circ \gamma_{j_1j_2} \circ \cdots \circ \gamma_{j_{k-1}j}(s_j) \\
	&> \omega \circ \omega^{-1}(\|s\|_{\ell_{\infty}}) = \|s\|_{\ell_{\infty}}.%
\end{align*}
The rest of the proof is the same as that of \cite[Lem.~B.5]{MKG20}.
\end{proof}

Some further simple properties of the operator $Q$ are summarized in the following proposition.%

\begin{proposition}\label{prop:Properties-of-Q}
Assume that $\Gamma:\ell_{\infty}^+ \rightarrow \ell_{\infty}^+$ is well-defined, continuous and satisfies the max-robust SGC. Then the operator $Q$ defined in \eqref{eq:Q-operator} has the following properties:
\begin{enumerate}[label = (\roman*)]
\item\label{itm:Simple-Q-properties-1} $Q(0) = 0$ and $Q$ is a monotone operator.%
\item\label{itm:Simple-Q-properties-2} The image of $Q$ is the set of all points of decay for $\Gamma$:%
\begin{equation*}
  \im\, Q = \{ s \in \ell_{\infty}^+ : \Gamma(s) \leq s \}.%
\end{equation*}
This set is closed, contains $s=0$, is cofinal (i.e., for any $x\in\ell_\infty^+$ there is $s \in\im\, Q$ with $x \leq s$) and forward-invariant with respect to $\Gamma$, i.e., $\Gamma(\im\, Q)\subset \im\, Q$.%
\item\label{itm:Simple-Q-properties-3} $Q \circ Q = Q$.%
\end{enumerate}
\end{proposition}

\begin{proof}
\ref{itm:Simple-Q-properties-1}. This immediately follows from the corresponding properties of $\Gamma$ and the definition of $Q$.%

\ref{itm:Simple-Q-properties-2}. Lemma~\ref{lem_q_op} implies that $\im\, Q \subset \{s \in \ell_{\infty}^+ : \Gamma(s) \leq s\}$. Conversely, $\Gamma(s) \leq s$ implies $\Gamma^k(s) \leq s$ for all $k \geq 0$, and hence, $Q(s) = s$ implying $s \in \im\, Q$. Since $\Gamma$ is continuous, it follows that $\im\, Q$ is closed. Since for each $s \in \ell_\infty^+$ it holds that $s \leq Q(s) \in \im\, Q $, the set $\im\, Q$ is cofinal. Since for any $s \in \im\, Q$ we have $\Gamma(s) \leq s$, by monotonicity of $\Gamma$ it follows that $\Gamma(\Gamma(s)) \leq \Gamma(s)$, showing forward-invariance of $\im\, Q$.%

\ref{itm:Simple-Q-properties-3}. This follows immediately from the proof of \ref{itm:Simple-Q-properties-2}.%
\end{proof}

For the gain operator $\Gamma$ and any $\theta \in \Kinf$, we define the operator $\Gamma_\theta: \ell_\infty^+ \to \ell_\infty^+$ by%
\begin{equation}\label{eq:A-omega}
  \Gamma_{\theta}(s) := (\id + \theta) \circ \Gamma(s) \mbox{\quad for all\ } s \in \ell_\infty^+.%
\end{equation}
Here, we apply the function $\id + \theta$ componentwise, i.e.%
\begin{equation*}
  \Gamma_{\theta}(s) = \big(\sup_{j \in \N} (\id + \theta) \circ \gamma_{ij}(s_j)\big)_{i \in \N}.%
\end{equation*}
Hence, the operator $\Gamma_\theta$ is structurally the same as $\Gamma$, but with scaled gain functions.%

We close the section with a simple lemma on gain operators satisfying the \emph{max-robust strong SGC}.%

\begin{lemma}\label{lem:Increasing-of-RS-SCL-operators} 
Assume that $\Gamma$ satisfies the max-robust strong SGC with some $\rho,\omega\in\Kinf$. Then $\Gamma_{\rho}$ satisfies the max-robust SGC with the same $\omega\in\Kinf$. Furthermore, there is $\theta \in \Kinf$ such that $\Gamma_{\theta}$ also satisfies the max-robust strong SGC.
\end{lemma}

\begin{proof}
Let $\Gamma$ satisfy the max-robust strong SGC with certain $\omega\in\Kinf$ and $\rho\in\Kinf$. By \cite[Lem.~1.1.3]{Rue07}, we can find $\rho_1,\rho_2\in\Kinf$ such that $\id + \rho = (\id + \rho_1) \circ (\id + \rho_2)$. Note that for any $s \in\ell_\infty^+$ and $i,j\in\N$ it holds that%
\begin{align*}
(\id +\rho)\circ \Gamma_{ij}(s) 
&= (\id+\rho_1) \circ (\id+\rho_2) \circ (\Gamma(s) \oplus \omega(s_j)e_i)\\
&\geq (\id+\rho_1) \circ ((\id+\rho_2)\circ \Gamma(s) \oplus \omega(s_j)e_i).
\end{align*}
Since $(\id +\rho) \circ \Gamma_{ij}(s) \not\geq s$ for all $s \neq 0$, it also holds that%
\begin{equation*}
  (\id + \rho_1) \circ (\Gamma_{\rho_2}(s) \oplus \omega(s_j)e_i)(s) \not\geq s,\quad \forall s \in \ell_{\infty}^+ \setminus \{0\},%
\end{equation*}
showing that $\Gamma_{\rho_2}$ satisfies the max-robust strong SGC. Similar arguments show that $\Gamma_{\rho}$ satisfies the max-robust SGC.
\end{proof}

\section{The discrete-time system induced by the gain operator}\label{sec:Discrete-time systems induced by the gain operator}

In this section, we relate the properties of the gain operator $\Gamma$ and its strong transitive closure $Q$ to the stability properties of the discrete-time system induced by the gain operator $\Gamma$:%
\begin{equation}
\label{eq:Gammaotimes-discrete-time-nonlinear}
  s(k+1) = \Gamma(s(k)),\quad k\in\Z_+.%
\end{equation}

As we will see, the stability properties of system \eqref{eq:Gammaotimes-discrete-time-nonlinear} play an important role in the construction of paths of strict decay for $\Gamma$, and hence for the construction of ISS Lyapunov functions for interconnected systems via the small-gain approach.%

The next proposition characterizes the max-robust SGC in terms of the stability properties of the system \eqref{eq:Gammaotimes-discrete-time-nonlinear}.%

\begin{proposition}[Criterion for max-robust SGC]
\label{prop_maxrobustsgc_char}
Let $\Gamma$ be well-defined and continuous. Then $\Gamma$ satisfies the max-robust SGC if and only if the following two properties hold:
\begin{enumerate}
\item[(i)] The system \eqref{eq:Gammaotimes-discrete-time-nonlinear} is \emph{uniformly globally stable (UGS)}, i.e., there is $\sigma\in\Kinf$ such that for any initial state $s\in \ell_\infty^+$, the solution of \eqref{eq:Gammaotimes-discrete-time-nonlinear} satisfies%
\begin{equation}\label{eq:UGS}
  \|\Gamma^k(s)\|_{\ell_\infty} \leq \sigma(\|s\|_{\ell_\infty}),\quad \forall k\in\Z_+.%
\end{equation}
\item[(ii)] Each trajectory of \eqref{eq:Gammaotimes-discrete-time-nonlinear} converges to zero componentwise, i.e.~$\pi_i \circ \Gamma^k(s) \rightarrow 0$ as $k \rightarrow \infty$ for every $s \in \ell_{\infty}^+$ and $i \in \N$.
\end{enumerate}
\end{proposition}

\begin{proof}
\q{$\Leftarrow$}: We first show that our assumption implies the existence of $\varphi \in \KC_{\infty}$ satisfying%
\begin{equation}\label{eq_intermediate_prop}
  s \leq \Gamma(s) \oplus b \quad \Rightarrow \quad \|s\|_{\ell_{\infty}} \leq \varphi(\|b\|_{\ell_{\infty}})%
\end{equation}
for any $s,b \in \ell_{\infty}^+$. To this end, assume that $s \leq \Gamma(s) \oplus b$. Using that $\Gamma$ is a max-preserving operator, inductively we obtain%
\begin{equation}\label{eq_sineq}
  s \leq \Gamma^k(s) \oplus \bigoplus_{l=0}^{k-1}\Gamma^l(b) \mbox{\quad for all\ } k \geq 1.%
\end{equation}
Looking at this inequality componentwise and letting $k \rightarrow \infty$ yields $s \leq \oplus_{l=0}^{\infty}\Gamma^l(b)$. Consequently, by the assumption of uniform global stability, there is $\varphi \in \Kinf$ such that%
\begin{equation*}
  \|s\|_{\ell_{\infty}} \leq \sup_{l \in \Z_+} \|\Gamma^l(b)\|_{\ell_{\infty}} \leq \varphi(\|b\|_{\ell_{\infty}}).%
\end{equation*}
Hence, the implication \eqref{eq_intermediate_prop} holds. To complete the proof, assume that $\Gamma(s) \oplus \omega(s_j)e_i \geq s$ for some $i,j \in \N$, $s \in \ell_{\infty}^+ \setminus \{0\}$ and $\omega < \varphi^{-1}$. Then $\|s\|_{\ell_{\infty}} \leq \varphi(\omega(s_j)) < s_j \leq \|s\|_{\ell_{\infty}}$, a contradiction.%

\q{$\Rightarrow$} 
By definition of $Q$, we have $\Gamma^k(s) \leq Q(s)$ for all $s \in\ell_\infty^+$, and by monotonicity of the norm and using Lemma~\ref{lem_q_op}, we obtain that $\|\Gamma^k(s)\|_{\ell_\infty} \leq \|Q(s)\|_{\ell_\infty} \leq \omega^{-1}(\|s\|_{\ell_\infty})$ for all $s,k$. This shows the UGS property.%

To show (ii), consider the operator $Q$ induced by $\Gamma$ and let $s \in \im\, Q$. Then $\Gamma(s) \leq s$ implying $\Gamma^{k+1}(s) \leq \Gamma^k(s)$ for all $k \in \Z_+$. Hence, each of the sequences $(\Gamma^k_i(s))_{k\in\Z_+}$, $i\in\N$, is monotonically decreasing and bounded below by zero. So the following limits exist:%
\begin{equation*}
  s^*_i := \lim_{k \rightarrow \infty}\Gamma^k_i(s),\quad i \in \N.%
\end{equation*}
The vector $s^* := (s^*_i)_{i\in\N}$ is an element of $\ell_{\infty}^+$, since $0 \leq s^* \leq s$. We claim that $s^* = 0$. To prove this, let $I_i$ be the finite set of $j \in \N$ with $\gamma_{ij} \neq 0$. Then observe that for each $i \in \N$%
\begin{align*}
  \Gamma_i(s^*) &= \sup_{j \in \N} \gamma_{ij}(s^*_j) = \sup_{j \in \N} \gamma_{ij}(\lim_{k\rightarrow\infty}\Gamma^k_j(s)) \\
								&\stackrel{(a)}{=} \sup_{j \in \N} \lim_{k \rightarrow \infty} \gamma_{ij}(\Gamma^k_j(s)) = \max_{j \in I_i} \lim_{k \rightarrow \infty} \gamma_{ij}(\Gamma^k_j(s)) \\
								&\stackrel{(b)}{=} \lim_{k \rightarrow \infty} \max_{j \in I_i} \gamma_{ij}(\Gamma^k_j(s)) = \lim_{k \rightarrow \infty} \Gamma_i(\Gamma^k(s)) \\
								&= \lim_{k \rightarrow \infty} \Gamma_i^{k+1}(s) = s^*_i.%
\end{align*}
The identity (a) holds because $\gamma_{ij}$ is continuous and (b) holds because the maximum (over finitely many quantities) commutes with the limit operation. Hence, the pointwise limit $s^*$ of the trajectory $(\Gamma^k(s))_{k\in\Z_+}$ is a fixed point of $\Gamma$. Since $\Gamma$ satisfies the SGC, this implies $s^* = 0$. We thus completed the proof for the case that $s \in \im\, Q$. For any other $s$, we have $s \leq Q(s)$, and hence $0 \leq \Gamma^k(s) \leq \Gamma^k(Q(s))$ for all $k \in \Z_+$, implying that $\Gamma^k(s) \rightarrow 0$ componentwise.
\end{proof}

In the next proposition, we define several basic stability properties of \eqref{eq:Gammaotimes-discrete-time-nonlinear} and show their equivalence.%

\begin{proposition}\label{prop:UGAS-criterion} 
Assume that $\Gamma$ satisfies the max-robust SGC. The following statements are equivalent:%
\begin{enumerate}[label = (\roman*)]
\item\label{itm:UGAS-criterion-sup-preserving-1} System \eqref{eq:Gammaotimes-discrete-time-nonlinear} is \emph{uniformly globally asymptotically stable (UGAS)}, i.e., there is $\beta\in\KL$, such that for any initial condition $s\in \ell_\infty^+$, the solution of \eqref{eq:Gammaotimes-discrete-time-nonlinear} satisfies%
\begin{equation}\label{eq:UGAS}
  \|\Gamma^k(s)\|_{\ell_\infty} \leq \beta(\|s\|_{\ell_\infty},k),\quad \forall k\in\Z_+.%
\end{equation}
\item\label{itm:UGAS-criterion-sup-preserving-2} System \eqref{eq:Gammaotimes-discrete-time-nonlinear} is \emph{globally attractive}, i.e., for all $s \in\ell_\infty^+$ it holds that $\Gamma^k(s)\to 0$ as $k\to\infty$.%
\item\label{itm:UGAS-criterion-sup-preserving-3} System \eqref{eq:Gammaotimes-discrete-time-nonlinear} is \emph{globally weakly attractive on $\im\, Q$}, i.e., $\inf_{k\geq 0}\|\Gamma^k(s)\|_{\ell_\infty} = 0$ for all $s \in \im\, Q$.
\end{enumerate}
\end{proposition}

\begin{proof}
Clearly, \ref{itm:UGAS-criterion-sup-preserving-1} $\Rightarrow$ \ref{itm:UGAS-criterion-sup-preserving-2} $\Rightarrow$ \ref{itm:UGAS-criterion-sup-preserving-3} holds. Let us show the implication \ref{itm:UGAS-criterion-sup-preserving-3} $\Rightarrow$ \ref{itm:UGAS-criterion-sup-preserving-1}.%

Let \eqref{eq:Gammaotimes-discrete-time-nonlinear} be globally weakly attractive on $\im\, Q$. For any $r>0$ and any $s \in B_r(0)$, it holds that $s \leq r \unit\le Q(r\unit) \in \im\, Q$. By monotonicity of $\Gamma$, it holds that $\Gamma^k(s) \leq \Gamma^k(Q(r \unit))$, and thus $\|\Gamma^k(s)\|_{\ell_\infty} \leq \|\Gamma^k(Q(r\unit))\|_{\ell_\infty}$ for all $k\in\Z_+$. Hence, $\inf_{k\geq 0}\sup_{s \in B_r(0)}\|\Gamma^k(s)\|_{\ell_\infty} = 0$ (the so-called uniform global weak attractivity of \eqref{eq:Gammaotimes-discrete-time-nonlinear}). Together with the UGS property, this implies UGAS (for continuous-time systems one can find this result, e.g., in \cite[Thm.~4.2]{Mir17a}; the proof of the discrete-time version is completely analogous).
\end{proof}

The next example shows that even if the gains $\gamma_{ij}$ are all linear and $\Gamma$ satisfies the strong as well as the max-robust SGC,
 the system induced by $\Gamma$ is not necessarily UGAS.

\begin{example}\label{sec:motivating-examples6}
Consider linear gains defined by%
\begin{equation*}
  \gamma_{k+1,k} := \delta_k \frac{k}{k+1} = 
\begin{cases}
0 & \text{ if } k\in \{2^s:s\in\N\},\\ 
\frac{k}{k+1} & \text{ otherwise.} 
\end{cases}
\end{equation*}
and $\gamma_{ij} = 0$ whenever $i \neq j+1$. Since $\gamma_{ij} \leq 1$ for all $k$, the gain operator $\Gamma$ is well-defined. Further, we have% 
\begin{align*}
  (\Gamma^{2^{k-1}-1}(\unit))_{2^k} &= 
  \gamma_{2^k,2^{k}-1} \cdot\ldots\cdot \gamma_{2^{k-1}+2,2^{k-1}+1} \\
  &= \frac{2^{k-1}+1}{2^k} \geq \frac{1}{2},%
\end{align*}
showing that $\Gamma^k(\unit)\not\to 0$ as $k\to\infty$, and thus the discrete-time system \eqref{eq:Gammaotimes-discrete-time-nonlinear} induced by $\Gamma$ is not UGAS. Now, \cite[Prop.~8]{homogeneous} implies, in particular, that \emph{$\Gamma$ does not satisfy the so-called robust strong SGC with linear $\rho$ and $\omega$}. On the other hand,%
\begin{itemize}
\item \emph{$\Gamma$ satisfies the strong SGC.} Assume to the contrary that there is $s\in\ell_\infty^+\setminus\{0\}$ such that $\Gamma s \geq (1-\ep)s$ for a fixed but arbitrary $\ep\in(0,1)$. Then we get $Rs \geq \Gamma s \geq (1-\ep)s$, where $R$ is the right-shift operator from \cite[Ex.~3.15]{GlM21}. However, in view of \cite[Ex.~3.15]{GlM21}, $R$ satisfies the strong SGC with any $\ep\in(0,1)$, and we come to a contradiction.%
\item \emph{$\Gamma$ satisfies the max-robust SGC.} Pick any $i,j \in \N$ and perturb the $ij$-component of $\Gamma$ by the linear function $\omega(r) = \frac{1}{2}r$. The resulting operator $\Gamma_{ij}$ is a block-diagonal operator with finite-dimensional blocks, and all its finite cycles are contractions. Thus, by \cite[Prop.~B.4]{MKG20}, $\Gamma_{ij}$ satisfies the SGC (here we use the fact that the robust SGC used in \cite{MKG20} implies the max-robust SGC). \hfill$\square$
\end{itemize}
\end{example}

We thus arrive at the following relations:%
\begin{figure}[h]
\begin{tikzpicture}[>=implies,thick]
\node (UGAS) at (-3.1,0.4) {UGAS};
\node (max-rob) at (-2.9,-0.6) {max-robust SGC};

\node (UGS) at (-7.7,0.4) {UGS $\wedge$ global attractivity (GATT)};
\node (UGS-w) at (-7.7,-0.6) {UGS $\wedge$ componentwise GATT};
\draw[thick,double equal sign distance,<->] (UGS) to (UGAS);
\draw[thick,double equal sign distance,<->] (UGS-w) to (max-rob);
\draw[thick,double equal sign distance,->] (-7,0.2) to (-7,-0.4);
\draw[thick,double equal sign distance,<-] (-6.3,0.2) to (-6.3,-0.4);
\node[thick,rotate=-90] (NEG) at (-6.3,-0.2) {$/$};
\end{tikzpicture}
\end{figure}
\hspace{-0.2cm}

Since we need UGAS for the construction of a path of strict decay, we need to understand what is required in addition to the max-robust SGC to obtain UGAS.%

It is well-known that for finite networks the max-preserving gain operator $\Gamma$ induces a UGAS discrete-time system if and only if all cycles composed of gains are contractions, see, e.g., \cite[Thm.~6.4]{Rue10}. In the case of infinite networks, UGAS of the induced system can be characterized in terms of sufficiently long chains of gains, as shown in the next proposition.

\begin{proposition}\label{prop_ugas_char}
Assume that the gain operator $\Gamma$ is well-defined, continuous and satisfies the max-robust small-gain condition. Then the following are equivalent:%
\begin{enumerate}
\item[(i)] The induced system \eqref{eq:Gammaotimes-discrete-time-nonlinear} is UGAS.%
\item[(ii)] There exists $\eta \in \KC$ with $\eta < \id$ such that for every $r \geq 0$ there is $n\in\N$ with%
\begin{equation}\label{eq_gaincontr}
  \sup_{j_0,j_1,\ldots,j_n\in\N} \gamma_{j_0j_1} \circ \cdots \circ \gamma_{j_{n-1}j_n}(r) \leq \eta(r).%
\end{equation}
\item[(iii)] There exist $\eta \in \KC$ with $\eta < \id$ and $i_0 \in \N$ such that for every $r \geq 0$ there is $n\in\N$ with%
\begin{equation*}
  \sup_{j_0,j_1,\ldots,j_n\in\N \atop j_0 \geq i_0} \gamma_{j_0j_1} \circ \cdots \circ \gamma_{j_{n-1}j_n}(r) \leq \eta(r).%
\end{equation*}
\end{enumerate}
\end{proposition}

\begin{proof}
(i) $\Rightarrow$ (ii): Assume that \eqref{eq:Gammaotimes-discrete-time-nonlinear} is UGAS. Then there exists $\beta \in \KC\LC$ such that $\|\Gamma^k(s)\|_{\ell_{\infty}} \leq \beta(\|s\|_{\ell_{\infty}},k)$ for all $s\in\ell_{\infty}^+$, $k \in \Z_+$. We put $\eta(r) := r/2$ for all $r \in \R_+$. For a given $r>0$, choose $n$ so large that $\beta(r,n) \leq r/2$. This implies%
\begin{align*}
 & \sup_{j_0,j_1,\ldots,j_n\in\N}\gamma_{j_0j_1} \circ \cdots \circ \gamma_{j_{n-1}j_n}(r) = \sup_{j_0\in\N} \Gamma^n_{j_0}(r\unit) \\
	&\qquad = \|\Gamma^n(r\unit)\|_{\ell_{\infty}} \leq \beta(r,n) \leq \eta(r).% 
\end{align*}
Hence, \eqref{eq_gaincontr} holds.%

(ii) $\Rightarrow$ (iii): This is obvious.%

(iii) $\Rightarrow$ (i): By Proposition \ref{prop:UGAS-criterion}, it suffices to prove global attractivity on $\im\, Q$. To this end, fix $s \in \im\, Q$ and put $r_1 := \|s\|_{\ell_{\infty}}$. By assumption, there exists $\tilde{n}_1 = \tilde{n}_1(r_1) \in \N$ with%
\begin{equation*}
  \sup_{j_0,j_1,\ldots,j_{\tilde{n}_1}\in\N \atop j_0 \geq i_0} \gamma_{j_0j_1} \circ \cdots \circ \gamma_{j_{\tilde{n}_1-1}j_{\tilde{n}_1}}(r_1) \leq \eta(r_1).%
\end{equation*}
This implies%
\begin{equation*}
  \Gamma^{\tilde{n}_1}_i(s) \leq \Gamma^{\tilde{n}_1}_i(r_1 \unit) \leq \eta(r_1) \mbox{\quad for all\ } i \geq i_0.%
\end{equation*}
By Proposition \ref{prop_maxrobustsgc_char}, we further find $\hat{n}_1 \in \N$ with%
\begin{equation*}
  \Gamma^{\hat{n}_1}_i(s) \leq \eta(r_1) \mbox{\quad for\ } 1 \leq i < i_0.%
\end{equation*}
Now, put $n_1 := \max\{\tilde{n}_1,\hat{n}_1\}$. 
As $s \in \im\, Q$, we obtain that $\Gamma^{n_1}(s) \leq \Gamma^{\hat{n}_1}(s)$ and $\Gamma^{n_1}(s) \leq \Gamma^{\tilde{n}_1}(s)$. Thus:
\begin{align*}
  \|\Gamma^{n_1}(s)\|_{\ell_{\infty}} &= \sup_{i\in\N}\Gamma^{n_1}_i(s) = \max\{ \max_{1 \leq i < i_0} \Gamma^{n_1}_i(s), \sup_{i \geq i_0} \Gamma^{n_1}_i(s) \} \\
	&\leq \max\{ \max_{1 \leq i < i_0} \Gamma_i^{\hat{n}_1}(s), \sup_{i \geq i_0} \Gamma_i^{\tilde{n}_1}(s) \} \leq \eta(r_1).%
\end{align*}
 Now, we put $r_2 := \eta(r_1)$. Proceeding in the same way, we find $n_2 \in \N$ with%
\begin{equation*}
  \|\Gamma^{n_1+n_2}(s)\|_{\ell_{\infty}} \leq \eta(r_2) = \eta^2(r_1).%
\end{equation*}
Inductively, we find a sequence $N_k \rightarrow \infty$ such that%
\begin{equation*}
  \|\Gamma^{N_k}(s)\|_{\ell_{\infty}} \leq \eta^k(r_1) \mbox{\quad for all\ } k \geq 0.%
\end{equation*}
Since $\eta < \id$ and $\eta \in \KC$, it follows that $\eta^k(r_1) \rightarrow 0$ (the only fixed point of $\eta$) as $k \rightarrow \infty$. This shows that \ref{itm:UGAS-criterion-sup-preserving-3} in Proposition \ref{prop:UGAS-criterion} is satisfied.%
\end{proof}

The UGAS property of \eqref{eq:Gammaotimes-discrete-time-nonlinear} implies important approximation and continuity properties of the operator $Q$.%

\begin{proposition}\label{prop:From-attractivity-to-Properties-of-Q-}
Assume that $\Gamma:\ell_{\infty}^+ \rightarrow \ell_{\infty}^+$ is well-defined and continuous and system \eqref{eq:Gammaotimes-discrete-time-nonlinear} is UGAS. Then the following statements hold:%
\begin{enumerate}[label = (\roman*)]
\item\label{itm:Attractivity-Q-properties-1} $\im\, Q$ is path-connected.%
\item\label{itm:Attractivity-Q-properties-2} For all $s_1,s_2 \in \inner(\ell_{\infty}^+)$ with $s_1\leq s_2$, there is $m\in\N$ such that%
\begin{equation}\label{eq:Q-representation-in-interior-of-the-cone}
  Q(s) = \bigoplus_{k = 0}^m \Gamma^k(s), \quad \forall s: s_1\leq s \leq s_2.%
\end{equation}
\item\label{itm:Attractivity-Q-properties-3} For each $q \in \inner(\ell_{\infty}^+)$ and $\varepsilon>0$, there is some $m\in\N$ such that for all $s \in\ell^+_\infty$ with $0\leq s \leq q$%
\begin{equation}\label{eq:Q-estimate-in-interior-of-the-cone-with-eps}
  \bigoplus_{k=0}^m\Gamma^k(s) \leq Q(s) \leq \bigoplus_{k=0}^m\Gamma^k(s) + \varepsilon\unit.%
\end{equation}
\item\label{itm:Attractivity-Q-properties-4} $Q$ is continuous on $\ell_{\infty}^+$.%
\end{enumerate}
\end{proposition}

\begin{proof}
\ref{itm:Attractivity-Q-properties-1}. It suffices to prove that any point $s \in \im\, Q$ can be connected with $0$ by a continuous path. Hence, take any $s \in \ell_{\infty}^+$ with $\Gamma(s) \leq s$. Now define $s_{\alpha} := \alpha s + (1 - \alpha)\Gamma(s)$ for all $\alpha \in [0,1]$ and observe that%
\begin{align*}
  \Gamma(s) &= \alpha \Gamma(s) + (1 - \alpha)\Gamma(s) \leq s_{\alpha} \\
	&= \alpha s + (1 - \alpha)\Gamma(s) \leq \alpha s + (1 - \alpha)s = s.%
\end{align*}
Thus, $\Gamma(s) \leq s_{\alpha} \leq s$ and applying $\Gamma$ once again and using its monotonicity yields $\Gamma(s_{\alpha}) \leq \Gamma(s) \leq s_{\alpha}$. This implies $s_{\alpha} \in \im\, Q$ for all $\alpha\in[0,1]$. We can repeat this process and thus construct a continuous (piecewise linear) path connecting $s$ with $\Gamma^k(s)$ for any $k \geq 1$. Since $\Gamma^k(s) \rightarrow 0$ by assumption, it follows that $s$ can be connected with $0$ by a continuous path.%

\ref{itm:Attractivity-Q-properties-2}. Pick any $s_1,s_2 \in \inner(\ell_{\infty}^+)$ with $s_1\leq s_2$. As $0$ is globally attractive, and $s_1 \in\inner(\ell_{\infty}^+)$, there is $m\in\N$ such that $\Gamma^k(s_2)\leq s_1$ for all $k\geq m$. Since $\Gamma$ is monotone, one gets $\Gamma^k(s)\leq \Gamma^k(s_2)\leq s_1 \leq s$ whenever $k\geq m$ and $s_1 \leq s \leq s_2$. This implies \eqref{eq:Q-representation-in-interior-of-the-cone}.%

\ref{itm:Attractivity-Q-properties-3}. Pick any $q \in \inner(\ell_{\infty}^+)$ and $\varepsilon>0$. As $0$ is globally attractive, there is $m\in\N$ such that $\Gamma^k(q) \leq \varepsilon \unit$ for all $k\geq m$. Since $\Gamma$ is monotone, then also $\Gamma^k(s)\leq \Gamma^k(q)\leq \varepsilon \unit$ whenever $k\geq m$ and $s \leq q$. By definition of $Q$, it holds that $Q(s) \leq \sup\{\bigoplus_{k=0}^m\Gamma^k(s), \varepsilon\unit\}$, which implies \eqref{eq:Q-estimate-in-interior-of-the-cone-with-eps}.%

\ref{itm:Attractivity-Q-properties-4}. Take any $q \in \inner(\ell_\infty^+)$, and pick $s\leq q$ satisfying $B_\delta(s) \subset [0,q] = \{x \in \ell_\infty^+:0\leq x\le q\}$ for some $\delta>0$. Take any $\tilde{s} \in B_{\delta}(s)$. By (iii), for each $\varepsilon>0$, we find $m\in\N$ such that%
\begin{equation*}
  \bigoplus_{k=0}^m\Gamma^k(w) \leq Q(w) \leq \bigoplus_{k=0}^m\Gamma^k(w) + \frac{1}{2}\varepsilon\unit,\quad \forall w \in [0,q],%
\end{equation*}
and thus%
\begin{align*}
  &\bigoplus_{k=0}^m\Gamma^k(s) - \bigoplus_{k=0}^m\Gamma^k(\tilde{s}) - \frac{1}{2}\varepsilon\unit \leq Q(s) - Q(\tilde{s}) \\
	&\qquad \leq \bigoplus_{k=0}^m\Gamma^k(s) - \bigoplus_{k=0}^m\Gamma^k(\tilde{s}) + \frac{1}{2}\varepsilon\unit.%
\end{align*}
Hence, 
\begin{equation*}
  \|Q(s) - Q(\tilde{s})\|_{\ell_\infty} \leq \Bigl\|\bigoplus_{k=0}^m\Gamma^k(s) - \bigoplus_{k=0}^m\Gamma^k(\tilde{s})\Bigr\|_{\ell_\infty} + \frac{\varepsilon}{2}.%
\end{equation*}
Using the simple representation of the iterates $\Gamma^k$ in terms of gains (see, e.g., \cite[Lem.~12]{MKG20}), we have%
\begin{align*}
  \|Q(s) &- Q(\tilde{s})\|_{\ell_{\infty}} \leq \sup_{i\in\N}\Bigl|\max_{k=0}^m \sup_{j_1,\ldots,j_k} \gamma_{ij_1}\circ\cdots\circ \gamma_{j_{k-1}j_k}(s_{j_k}) \\
	& - \max_{k=0}^m \sup_{j_1,\ldots,j_k} \gamma_{ij_1}\circ\cdots\circ \gamma_{j_{k-1}j_k}(\tilde{s}_{j_k})\Bigr| + \frac{\varepsilon}{2}.%
\end{align*}
Since $\{\gamma_{ij}:i,j\in\N\}$ is pointwise equicontinuous by assumption, the same is true for the family $\{\gamma_{ij_1}\circ\cdots\circ\gamma_{j_{k-1}j_k} : i,j_1,\ldots,j_k \in \N,\ 0 \leq k < m\}$. Hence, this family is uniformly equicontinuous on compact intervals, and thus there exists $\rho \in (0,\delta)$ such that $\|s - \tilde{s}\|_{\ell_{\infty}} \leq \rho$ implies%
\begin{equation*}
  |\gamma_{ij_1}\circ\cdots\circ \gamma_{j_{k-1}j_k}(s_{j_k}) - \gamma_{ij_1}\circ\cdots\circ \gamma_{j_{k-1}j_k}(\tilde{s}_{j_k})| \leq \frac{1}{2}\ep%
\end{equation*}
for all $i,j_1,\ldots,j_k \in \N$ and $0 \leq k < m$. This, in turn, implies $\|Q(s) - Q(\tilde{s})\|_{\ell_{\infty}} \leq \ep$. We have thus proved continuity of $Q$ at $s$, and since any $s \in \ell_{\infty}^+$ satisfies $s \leq q$ for some $q \in \inner(\ell_{\infty}^+)$, the proof is complete.
\end{proof}

\section{Construction of paths of strict decay}\label{sec:Paths of strict decay}

We can finally present our main result on the existence and construction of paths of strict decay. It extends the first result of this kind in \cite[Lem.~4.5]{dashkovskiy2019stability}, where properties (i)--(iii) of a path of strict decay have been shown under similar assumptions. See Section~\ref{sec:Discussion of obtained results} for an extended discussion of this issue.%

%Then for each compact interval $K \subset (0,\infty)$, there exist $0 < c \leq C < \infty$ such that%
%\begin{equation*}
%  c|r_1 - r_2| \leq |\sigma_i^{-1}(r_1) - \sigma_i^{-1}(r_2)| \leq C|r_1 - r_2|%
%\end{equation*}
%for all $r_1,r_2 \in K$ and $i \in \N$.%

\begin{theorem}
\label{thm_omegapath}
Let the following assumptions hold:%
\begin{enumerate}
\item[(a)] There exists $\theta \in \Kinf$ such that the system induced by $\Gamma_{\theta} = (\id + \theta) \circ \Gamma$ is UGAS.%
\item[(b)] For each compact interval $K \subset (0,\infty)$, there are $0 < l \leq L < \infty$ with $l(r_2 - r_1) \leq \gamma_{ij}(r_2) - \gamma_{ij}(r_1) \leq L(r_2 - r_1)$ for all nonzero $\gamma_{ij}$ and $r_1 < r_2$ in $K$. 
\end{enumerate}
Then there exists a path of strict decay $\sigma:\R_+ \rightarrow \ell_{\infty}^+$ for $\Gamma$.
\end{theorem}

\begin{proof}
First, we fix $\theta\in\Kinf$ such that the system induced by $\Gamma_{\theta}$ is UGAS. We also put $\gamma_{ij}^{\theta} := (\id + \theta)\circ\gamma_{ij}$ for all $i,j\in\N$, and define $\sigma(r) := Q_{\theta}(r\unit)$ for all $r \in \R_+$, where $Q_{\theta}(s) = \bigoplus_{k\in\Z_+}\Gamma_{\theta}^k(s)$. Then we can verify all properties of a path of strict decay for $\sigma$:%

Since $\Gamma_{\theta}(Q_{\theta}(r\unit)) \leq Q_{\theta}(r\unit)$ by Lemma \ref{lem_q_op}, we obtain $\Gamma(\sigma(r)) \leq (\id + \theta)^{-1} \circ \sigma(r)$ for all $r \in \R_+$. Hence, $\sigma$ satisfies property (i) of a path of strict decay.%

Property (ii) of a path of strict decay holds with $\sigma_{\min} = \id$, since $\sigma(r) \geq r\unit$, and $\sigma_{\max} = \omega^{-1}$ by Lemma \ref{lem_q_op}.%

From property (ii), we can conclude that $\sigma_i(0) = 0$, $\sigma_i(r) > 0$ for all $r > 0$, and $\sigma_i(r) \rightarrow \infty$ as $r \rightarrow \infty$. As $0$ is a globally attractive fixed point for $\Gamma_{\theta}$, $Q_{\theta}$ is continuous on $\ell_\infty^+$ by Proposition~\ref{prop:From-attractivity-to-Properties-of-Q-}. Thus, all $\sigma_i$ are continuous as well. Furthermore, for $r_1,r_2 \in(0,\infty)$ with $r_1<r_2$, we obtain by Proposition~\ref{prop:From-attractivity-to-Properties-of-Q-}\ref{itm:Attractivity-Q-properties-2}, that%
\begin{equation*}
  \sigma_i(r) = \max_{0 \leq k < k_0}\pi_i \circ \Gamma_{\theta}^k(r\unit) \mbox{\quad for all\ } r \in [r_1,r_2].%
\end{equation*}
By our assumption that for each $i$ only finitely many $\gamma_{ij}$ are nonzero, the supremum in%
\begin{equation*}
  \sigma_i(r) = \max_{0 \leq k < k_0}\sup_{j_1,\ldots,j_k} \gamma^{\theta}_{ij_1} \circ \cdots \circ \gamma^{\theta}_{j_{k-1}j_k}(r)%
\end{equation*}
is in fact a supremum over finitely many strictly increasing functions (since we can ignore all chains which contain a zero function). This implies that $\sigma_i$ is also strictly increasing on $[r_1,r_2]$, and hence everywhere. It follows that all $\sigma_i$ are $\Kinf$-functions (property (iii) of a path of strict decay).%

We complete the proof by verifying property (iv) of a path of strict decay. It suffices to prove the statement for $\sigma_i$ in place of $\sigma_i^{-1}$. Indeed, assume that for every compact interval $L \subset (0,\infty)$, we have constants $\tilde{c},\tilde{C}>0$ satisfying%
\begin{equation*}
  \tilde{c}|r_1 - r_2| \leq |\sigma_i(r_1) - \sigma_i(r_2)| \leq \tilde{C}|r_1 - r_2|\ \forall r_1,r_2 \in L.%
\end{equation*}
If $K = [a,b] \subset (0,\infty)$, then $\sigma_i^{-1}(K) = [\sigma_i^{-1}(a),\sigma_i^{-1}(b)]$ which is a subset of $[\sigma_{\max}^{-1}(a),\sigma_{\min}^{-1}(b)] =: L \subset (0,\infty)$. Hence, the above estimates imply%
\begin{equation*}
  \frac{1}{\tilde{C}}|r_1 - r_2| \leq |\sigma_i^{-1}(r_1) - \sigma_i^{-1}(r_2)| \leq \frac{1}{\tilde{c}}|r_1 - r_2|\ \forall r_1,r_2 \in K.%
\end{equation*}
To verify the statement for $\sigma_i$, we first prove the following claim: If $K = [a,b] \subset (0,\infty)$ is a compact interval, then there exists another compact interval $[c,d] \subset (0,\infty)$ such that $\gamma_{ij}^{\theta}(K) \subset [c,d]$ for all nonzero $\gamma_{ij}$. By uniform equicontinuity of the functions $\{\gamma_{ij}^{\theta}\}$ on compact intervals, we know that $\gamma_{ij}^{\theta}(b)$ is uniformly bounded from above by some $d>0$. Now pick $\rho>0$ such that $a-\rho > 0$. Then by assumption, there is $l>0$ with $\gamma_{ij}^{\theta}(a) - \gamma_{ij}^{\theta}(a-\rho) \geq l\rho$ for all $i,j \in \N$ such that $\gamma_{ij} \neq 0$. Hence, $\gamma_{ij}^{\theta}(a) \geq (1+\theta)l\rho$ whenever $\gamma_{ij} \neq 0$. This implies $\gamma_{ij}^{\theta}(K) \subset [c,d]$ with $c := (1+\theta)l\rho$.%

Now fix a compact interval $K \subset (0,\infty)$, $k \in \N$, and consider all chains of the form%
\begin{equation*}
  c_{j_1\ldots j_k} := \gamma^{\theta}_{j_1j_2} \circ \gamma^{\theta}_{j_2j_3} \circ \cdots \circ \gamma^{\theta}_{j_{k-1}j_k}%
\end{equation*}
which are built from nonzero gains. From the claim and our assumptions, it then follows that%
\begin{align*}
  l_1l_2 \cdots l_k |r_1 - r_2| &\leq |c_{j_1\ldots j_k}(r_1) - c_{j_1\ldots j_k}(r_2)| \\
	&\leq L_1L_2 \cdots L_k |r_1 - r_2|%
\end{align*}
for certain positive numbers $l_i,L_i>0$, $i = 1,\ldots,k$, and all $r_1,r_2 \in K$. The same Lipschitz bounds then also hold for the functions $r \mapsto \pi_i \circ \Gamma^k_{\theta}(r\unit) = \sup_{j_2\ldots j_k}c_{ij_2\ldots j_k}(r)$, where for the lower bound we need to require that at least one nonzero chain $c_{ij_2\ldots j_k}$ exists. By what we have shown above, on every compact interval $K \subset (0,\infty)$, $\sigma_i$ can be written as the maximum over finitely many of such functions:%
\begin{equation}\label{eq_sigmai_max}
  \sigma_i(r) = \max_{0 \leq k < k_0}\max_{j_2\ldots j_k}c_{ij_2\ldots j_k}(r) \mbox{\quad for all\ } r \in K.%
\end{equation}
With $C := \max\{1,L_1,L_1L_2,\ldots,L_1\cdots L_{k_0}\}$ (keeping in mind that $\pi_i \circ \Gamma_{\theta}^0(r\unit) = r$), then the following holds:%
\begin{equation*}
  |\sigma_i(r_1) - \sigma_i(r_2)| \leq C|r_1 - r_2| \mbox{\quad for all\ } r_1,r_2 \in K.%
\end{equation*}
For the lower bound, we put $c := \min\{1,l_1,l_1l_2,\ldots,l_1l_2\cdots l_{k_0}\}$ and apply Lemma \ref{lem_incrmax}. Observe that in taking the supremum, we can ignore all functions which are identically zero. Considering $k=0$, we see that at least one nonzero function is involved in taking the supremum, namely, the identity.%
\end{proof}

%Some comments regarding the assumptions of Theorem \ref{thm_omegapath} are necessary:%

%First, Lemma \ref{lem:Increasing-of-RS-SCL-operators} shows that the required robust small-gain condition for the operator %$\Gamma_{\theta}$ can be guaranteed if $\Gamma$ satisfies the robust strong small-gain condition and $\eta$ is chosen accordingly.%

\begin{remark}
Properties (i)--(iii) of a path of strict decay actually hold for the mapping $\sigma(r) = Q_{\theta}(r\unit)$ under the assumption that $\Gamma_{\theta}$ satisfies the max-robust SGC only, since this already implies that $\pi_i \circ \Gamma_{\theta}^k(r\unit) \rightarrow 0$ as $k \rightarrow \infty$ for every $i \in \N$, which is enough to write $\sigma_i$ locally as the maximum over finitely many continuous and strictly increasing functions. The stronger assumption of UGAS is only needed to verify the uniform local Lipschitz condition (iv), which requires that we can locally write $\sigma_i(r) = \max_{0 \leq k < k_0} \pi_i \circ \Gamma^k_{\theta}(r\unit)$ with $k_0$ being independent of $i$.
\end{remark}

\begin{remark}
Assumption (b) in Theorem \ref{thm_omegapath} is, in general, unnecessarily strong. This can be seen best by looking at the case when all gains $\gamma_{ij}$ are linear functions. Then the maximum in \eqref{eq_sigmai_max} is taken over linear functions, and since one of them is the identity, only those linear functions with slope $>1$ need to be taken into account. However, then it is clear that no uniform lower bound on the slopes of the individual $\gamma_{ij}$ are necessary to obtain a uniform lower bound on the Lipschitz constants of the $\sigma_i$. In general, the situation is more complicated, and detailed information about the gains and their compositions is necessary to relax assumption (b).
\end{remark}

\section{Sufficient conditions for UGAS}\label{sec:Sufficient conditions for attractivity}

Since UGAS of the discrete-time system, induced by the scaled gain operator, is a key requirement for the existence of a path of strict decay (and thus, for the application of the small-gain theorem), in this section, we analyze sufficient conditions for UGAS of the system induced by the gain operator $\Gamma$. The next proposition describes a way of reducing the proof of UGAS of \eqref{eq:Gammaotimes-discrete-time-nonlinear} to finitely many computations.%

\begin{proposition}
Assume that there exists a positive integer $N$ and a map $p:\N \rightarrow \{1,\ldots,N\}$ as well as a family $\{\bar{\gamma}_{ij} : i,j=1,\ldots,N\} \subset \KC \cup \{0\}$ of virtual gains such that%
\begin{equation*}
  \gamma_{ij} \leq \bar{\gamma}_{p(i)p(j)} \mbox{\quad for all\ } i,j \in \N.%
\end{equation*}
Let $\bar{\Gamma}:\R^N_+ \rightarrow \R^N_+$, $s \mapsto (\sup_{1\leq j \leq N}\bar{\gamma}_{ij}(s_j))_{1 \leq i \leq N}$ be the associated virtual gain operator. If $\bar{\Gamma}$ satisfies the SGC, $\bar{\Gamma}(s) \not\geq s$ for all $s \in \R^N_+ \setminus \{0\}$, the system \eqref{eq:Gammaotimes-discrete-time-nonlinear} induced by $\Gamma$ is UGAS.
\end{proposition}

\begin{proof}
For any $s \in \ell_{\infty}^+$ and $i \in \N$, we have with a convention $j_0:=i$:
\begin{align*}
  \Gamma^k_i(s) &= \sup_{j_1,\ldots,j_k \in \N} \gamma_{ij_1} \circ \cdots \circ \gamma_{j_{k-1}j_k}(s_{j_k}) \\
	&\leq \sup_{j_1,\ldots,j_k \in \N} \bar{\gamma}_{p(i)p(j_1)}\circ \cdots \circ \bar{\gamma}_{p(j_{k-1})p(j_k)}(\|s\|_{\ell_{\infty}}) \\
	&\leq \sup_{1 \leq j_1,\ldots,j_k \leq N} \bar{\gamma}_{p(i) j_1} \circ \bar{\gamma}_{j_2j_3} \circ \cdots \circ \bar{\gamma}_{j_{k-1}j_k}(\|s\|_{\ell_{\infty}}) \\
	&= \bar{\Gamma}^k_{p(i)}(\|s\|_{\ell_{\infty}}\unit),%
\end{align*}
where $\unit = (1,\ldots,1) \in \R^N_+$ in the last line. By \cite[Thm.~6.4]{Rue10}, the assumption that $\bar{\Gamma}$ satisfies the SGC implies that the $\ell_{\infty}$-norm in $\R^N$ of $\bar{\Gamma}^k(\|s\|_{\ell_{\infty}}\unit)$ is bounded by $\beta(\|s\|_{\ell_{\infty}},k)$ for some $\KL$-function $\beta$ (depending on $\bar{\Gamma}$ only). Hence, the system induced by $\Gamma$ is UGAS.
\end{proof}

Another method of checking UGAS of \eqref{eq:Gammaotimes-discrete-time-nonlinear} via the introduction of virtual gains, based on a compactification of the index set $\N$, is described in the next proposition.%

\begin{proposition}\label{prop_compactness_crit}
Let $\N^* := \N \cup \{\infty\}$ and assume that there exist virtual gains $\bar{\gamma}_{ij} \in \KC \cup \{0\}$, $i,j \in \N^*$ (where $\bar{\gamma}_{\infty\infty} \neq 0$ is allowed), satisfying the following assumptions:%
\begin{enumerate}
\item[(i)] $\bar{\gamma}_{ij} = \gamma_{ij}$ whenever $(i,j) \in \N \tm \N$.%
\item[(ii)] The virtual gain operator%
\begin{equation*}
  \bar{\Gamma}:\ell_{\infty}^+(\N^*) \rightarrow \ell_{\infty}^+(\N^*),\quad s \mapsto (\sup_{j\in\N^*} \bar{\gamma}_{ij}(s_j))_{i\in\N^*},%
\end{equation*}
is well-defined, continuous and satisfies the max-robust SGC with some $\omega \in \KC_{\infty}$.%
\item[(iii)] For each $i\in\N^*$, $\bar{\gamma}_{ij} \neq 0$ only for finitely many $j \in \N^*$.%
\item[(iv)] There exists $k_0 \in \N$ such that for all $r > 0$%
\begin{align}\label{eq_usc_cond}
\begin{split}
  &\limsup_{i \rightarrow \infty}\sup_{j_1,\ldots,j_{k_0} \in \N^*} \bar{\gamma}_{ij_1} \circ \cdots \circ \bar{\gamma}_{j_{k_0-1}j_{k_0}} \circ \omega^{-1}(r) \\
	&\leq \sup_{j_1,\ldots,j_{k_0} \in \N^*} \bar{\gamma}_{\infty j_1} \circ \cdots \circ \bar{\gamma}_{j_{k_0-1}j_{k_0}}(r).%
\end{split}
\end{align}
\end{enumerate}
Then, the system \eqref{eq:Gammaotimes-discrete-time-nonlinear} induced by $\Gamma$ is UGAS.%
\end{proposition}

\begin{proof}
First, we turn $\N^*$ into a compact metrizable space by declaring $\infty$ to be an accumulation point and all other elements isolated points. Now, let $s = \bar{Q}(r\unit)$ for some $r>0$, where $\bar{Q}(s) = \sup_{k\geq0}\bar{\Gamma}^k(s)$. Then, by Lemma \ref{lem_q_op}, $\bar{\Gamma}(s) \leq s$ and $r\unit \leq s \leq \omega^{-1}(r)\unit$. For every $k \in \Z_+$, we define%
\begin{equation*}
  f_k(i) := \bar{\Gamma}^k_i(s),\quad f_k:\N_* \rightarrow \R_+.%
\end{equation*}
By the choice of $s$, the functions $f_k$ form a decreasing sequence, i.e.~$f_{k+1} \leq f_k$ for all $k\in\Z_+$. We claim that the function $f_{k_0}$ is upper semicontinuous. Since $\infty$ is the only accumulation point of $\N_*$, this is equivalent to%
\begin{equation*}
  \limsup_{i \rightarrow \infty} \bar{\Gamma}^{k_0}_i(s) \leq \bar{\Gamma}^{k_0}_{\infty}(s).%
\end{equation*}
In terms of the gains $\bar{\gamma}_{ij}$, this can be written as%
\begin{align*}
  &\limsup_{i \rightarrow \infty} \sup_{j_1,\ldots,j_{k_0} \in \N^*} \bar{\gamma}_{ij_1} \circ \cdots \circ \bar{\gamma}_{j_{k_0-1}j_{k_0}}(s_{j_{k_0}}) \\
	&\leq \sup_{j_1,\ldots,j_{k_0} \in \N^*} \bar{\gamma}_{\infty j_1} \circ \cdots \circ \bar{\gamma}_{j_{k_0-1}j_{k_0}}(s_{j_{k_0}}).%
\end{align*}
From $r\unit \leq s \leq \omega^{-1}(r)\unit$, it follows that this inequality is implied by our assumption \eqref{eq_usc_cond}. Hence, $f_{k_0}$ is upper semicontinuous. It now easily follows that also $f_{kk_0}$ is upper semicontinuous for every $k \in \N$. Consequently, $(f_{kk_0})_{k\in\Z_+}$ is a decreasing sequence of upper semicontinuous functions on a compact metric space. Hence, Lemma \ref{lem_supinf} implies%
\begin{align*}
  & \lim_{k \rightarrow \infty} \|\bar{\Gamma}^{kk_0}(s)\|_{\ell_{\infty}} = \inf_{k\in\Z_+} \sup_{i \in \N^*} f_{kk_0}(i) \\
	&= \sup_{i \in \N^*} \inf_{k\in\Z_+} f_{kk_0}(i) = \sup_{i \in \N^*} \lim_{k \rightarrow \infty} \bar{\Gamma}^{kk_0}_i(s) = 0,%
\end{align*}
where the last equality follows from Proposition \ref{prop_maxrobustsgc_char}. This clearly implies $\bar{\Gamma}^k(s) \rightarrow 0$ for $k \rightarrow \infty$. Now, it is easy to see that the same property holds for the original gain operator $\Gamma$. The proof is completed by applying Proposition \ref{prop:UGAS-criterion}.%
\end{proof}

\begin{remark}
To get a better understanding of condition \eqref{eq_usc_cond}, consider the case $k_0=1$:%
\begin{align*}
  \limsup_{i \rightarrow \infty}\sup_{j \in \N^*} \bar{\gamma}_{ij} \circ \omega^{-1}(r) \leq \sup_{j \in \N^*} \bar{\gamma}_{\infty j}(r).%
\end{align*}
This implies that for all sufficiently large $i$, up to some small error, one gets $\bar{\gamma}_{ij} \circ \omega^{-1} \leq \bar{\gamma}_{\infty j^*}$ for some $j^*$ and all $j \in \N^*$. However, from the proof of Lemma \ref{lem_q_op}, we know that $\bar{\gamma}_{\infty j^*} < \omega^{-1}$ as a consequence of the max-robust SGC. Hence,%
\begin{equation*}
  \bar{\gamma}_{ij} \leq \bar{\gamma}_{\infty j^*} \circ \omega < \omega^{-1} \circ \omega = \id.%
\end{equation*}
Thus, asymptotically, the gains become contractions. If $k_0 > 1$, this is true for all chains of $k_0$ gains.
\end{remark}

\begin{example}
Consider a cascade network built from subsystems%
\begin{align*}
  \Sigma_1:&\quad \dot{x}_1 = f_1(x_1,u_1), \\
  \Sigma_i:&\quad \dot{x}_i = f_i(x_i,x_{i-1},u_i),\quad i \geq 2.%
\end{align*}
We assume that ISS Lyapunov functions $V_i$ for $\Sigma_i$ exist with associated interconnection gains $\gamma_{i(i-1)}$, which form an equicontinuous family. Then, the gain operator $\Gamma$ is well-defined and continuous. Let us assume that $\Gamma$ satisfies the max-robust small-gain condition with some $\omega \in \KC_{\infty}$, $\omega < \id$. We add a virtual gain $\bar{\gamma}_{\infty \infty} \in \KC \cup \{0\}$. It is then easy to see that with $\bar{\gamma}_{i(i-1)} := \gamma_{i(i-1)}$ for all $i\in\N$, $i \geq 2$, the virtual gain operator $\bar{\Gamma}$ satisfies the max-robust small-gain condition with the same $\omega$ if and only if%
\begin{equation*}
  \bar{\gamma}_{\infty\infty}(r) < r \mbox{\quad for all\ } r > 0.%
\end{equation*}
Additionally, the asymptotic condition (iv) in Proposition \ref{prop_compactness_crit} for $k_0 = 1$ requires that%
\begin{equation*}
  \limsup_{i \rightarrow \infty}\gamma_{i(i-1)} \circ \omega^{-1}(r) \leq \bar{\gamma}_{\infty \infty}(r) \mbox{\quad for all\ } r > 0.%
\end{equation*}
This means that the gains $\gamma_{i(i-1)}$ asymptotically have to become smaller than the composition of the two contractions $\bar{\gamma}_{\infty\infty}$ and $\omega$.%
\end{example}

The following proposition characterizes an interesting special case of Proposition \ref{prop_compactness_crit}.%

\begin{proposition}
Assume that the virtual gains in Proposition \ref{prop_compactness_crit} are chosen as $\bar{\gamma}_{ij} := 0$ if $i = \infty$ or $j = \infty$. Further, assume that $\Gamma$ is well-defined, continuous and satisfies the max-robust small-gain condition. Then the same is true for the virtual gain operator $\bar{\Gamma}$. Moreover, the following implications hold:%
\begin{enumerate}
\item[(i)] If \eqref{eq_usc_cond} holds, then the operator $\Gamma^{k_0}$ is compact, i.e., the image of any bounded set under $\Gamma^{k_0}$ is relatively compact.%
\item[(ii)] If the operator $\Gamma^{k_0}$ is compact and every subsystem of $\Sigma$ can only influence finitely many other subsystems, then condition \eqref{eq_usc_cond} is satisfied.%
\end{enumerate}
\end{proposition}

\begin{proof}
The proof is subdivided into three steps.%

\emph{Step 1}: The proof that $\bar{\Gamma}$ is well-defined or continuous if $\Gamma$ has the corresponding property is trivial, and hence we omit it. To verify that $\bar{\Gamma}$ satisfies the max-robust small-gain condition if $\Gamma$ does, we need to verify that for all $i,j \in \N^*$, the operator $\bar{\Gamma}_{ij}(s) = \bar{\Gamma}(s) \oplus \omega(s_j)e_i$ satisfies $\bar{\Gamma}_{ij}(s) \not\geq s$ for all $s \in \ell_{\infty}^+(\N^*) \setminus \{0\}$. To this end, we fix $s = (s_i)_{i\in\N^*}$ and distinguish several cases:%
\begin{itemize}
\item $i,j < \infty$ and $s_k > 0$ for some $k \in \N$. Then, with $\hat{s} := (s_i)_{i\in\N}$, the claim follows directly from $\Gamma_{ij}(\hat{s}) \not\geq \hat{s}$.%
\item $i,j < \infty$ and $s_k = 0$ for all $k \in \N$. Then $s_{\infty} > 0$ and thus $\bar{\Gamma}_{\infty}(s) = 0 < s_{\infty}$, showing that $\bar{\Gamma}_{ij}(s) \not\geq s$.%
\item $i = \infty$ or $j = \infty$ and $s_k > 0$ for some $k \in \N$. Then, with $\hat{s} := (s_i)_{i\in\N}$, it follows that $\Gamma(\hat{s}) \not\geq \hat{s}$, and consequently $\bar{\Gamma}_{ij}(s) \not\geq s$.%
\item ($i = \infty$, $j < \infty$) or ($i < \infty$, $j = \infty$) and $s_k = 0$ for all $k \in \N$. In both cases, $\pi_{\infty} \circ \bar{\Gamma}_{ij}(s) = 0 < s_{\infty}$.%
\item $i,j = \infty$ and $s_k = 0$ for all $k \in \N$. Then $\omega < \id$ implies%
\begin{equation*}
  \bar{\Gamma}_{\infty}(s) \oplus \omega(s_{\infty}) = \omega(s_{\infty}) < s_{\infty}.%
\end{equation*}
\end{itemize}

\emph{Step 2}: We prove statement (i). Without loss of generality, we assume that $k_0 = 1$. Condition \eqref{eq_usc_cond} can then be reformulated as follows: for each $r > 0$ and $\ep > 0$, there is $i_0 = i_0(r,\ep) \in \N$ such that for all $i \geq i_0$ and $j \in \N$ we have $\gamma_{ij}(r) \leq \ep$. Now, let $B \subset \ell_{\infty}^+$ be a bounded set, say $\|s\|_{\ell_{\infty}} \leq b$ for all $s \in B$. We have to show that $\Gamma(B)$ is relatively compact. By \cite[Thm.~6, p.~260]{dunford1988linear}, this is equivalent to the following: $\Gamma(B)$ is bounded and for every $\ep>0$, there is a finite partition $\N = J_1 \cup \ldots \cup J_N$ such that for every $s \in \Gamma(B)$ we have $|s_i - s_j| \leq \ep$ whenever $i,j \in J_k$ for some $k$. Boundedness of $\Gamma(B)$ is a consequence of the equicontinuity of $\{\gamma_{ij}\}$. If we choose $i_0 = i_0(b,\ep)$ as above, we can put $J_l := \{l\}$ for $l = 1,\ldots,i_0-1$ and $J_{i_0} := \{ i \in \N : i \geq i_0 \}$. For any $s = \Gamma(t) \in \Gamma(B)$, we have%
\begin{align*}
  s_{i_1} - s_{i_2} &= \sup_{j \in \N} \gamma_{i_1 j}(t_j) - \sup_{j \in \N} \gamma_{i_2 j}(t_j) \leq \sup_{j \in \N} \gamma_{i_1 j}(t_j) \\
	&\leq \sup_{j \in \N} \gamma_{i_1 j}(b) \leq \ep%
\end{align*}
whenever $i_1 \geq i_0$. Interchanging the roles of $i_1$ and $i_2$ yields $|s_{i_1} - s_{i_2}| \leq \ep$ whenever $i_1,i_2 \in J_{i_0}$. Hence, $\Gamma$ is compact.%

\emph{Step 3}: We prove statement (ii) for $k_0 = 1$. To this end, let us assume that condition \eqref{eq_usc_cond} does not hold. Then there exists $r>0$ such that%
\begin{equation*}
  \limsup_{i \rightarrow \infty} \sup_{j \in \N} \gamma_{ij}(r) > 0.%
\end{equation*}
Hence, there exist $K_0\in\N$, $\alpha > 0$ and sequences $i_k \rightarrow \infty$, $j_k \in I_{i_k}$ such that $\gamma_{i_kj_k}(r) \geq \alpha$ for all $k \geq K_0$. Now, consider the bounded set $C := \{r e_{j_k} : k \in \N\} \subset \ell_{\infty}^+$. We prove that $\Gamma(C)$ is not relatively compact, and thus $\Gamma$ is not compact. First, observe that%
\begin{equation*}
  \Gamma_{i_k}(r e_{j_k}) \geq \gamma_{i_k j_k}(r) \geq \alpha \mbox{\quad for all\ } k \geq K_0.%
\end{equation*}
Moreover, by our assumption that every subsystem can only influence finitely many other subsystems, we have $\Gamma_i(r e_{j_k}) = 0$ for almost all $i \in \N$. Assume to the contrary that $\N = J_1 \cup \ldots \cup J_N$ is a finite partition such that $|\Gamma_i(r e_{j_k}) - \Gamma_j(r e_{j_k})| \leq \alpha/2$ whenever $i,j \in J_l$ for some $l \in \{1,\ldots,N\}$. Assuming that $J_l$ contains infinitely many elements (which must be true for at least one of the index sets), we find that $\Gamma_{i_k}(re_{j_k}) - \Gamma_j(re_{j_k}) \geq \alpha - 0 = \alpha$ for some $j\in J_l$, a contradiction.% 
\end{proof}

\section{Discussion of the obtained results}\label{sec:Discussion of obtained results}

\textbf{Comparison to \cite{DaP19}.} The first Lyapunov-based small-gain theorem for infinite networks can be found in \cite[Thm.~1]{DaP19}. It requires that all internal gains are identical ($\gamma_{ij} \equiv \gamma$), and the SGC in \cite{DaP19} requires that $\gamma < \id$, which is in most cases much more conservative than the SGC employed in our paper.%

As $\gamma_{ij}=\gamma<\id$ for all $i,j\in\N$, the gain operator in \cite{DaP19} has a simple representation%
\begin{equation*}
  \Gamma(s) = \Bigl( \sup_{j\in\N\setminus\{i\}}\gamma(s_j) \Bigr)_{i\in\N} \leq \gamma(\|s\|_{\ell_\infty})\unit.%
\end{equation*}
Define $\sigma:\R_+\to\ell_\infty^+$ by $\sigma(r):=r\unit$, $r\geq0$. Then%
\begin{equation}\label{eq:Non-uniform decay}
  \Gamma(\sigma(r)) \leq \gamma(\|r\unit\|_{\ell_\infty})\unit = \gamma(r)\unit = \gamma\circ \sigma(r) < \sigma(r).%
\end{equation}
If $\gamma<\id$ satisfies%
\begin{equation}\label{eq:Paths-of-strict-decay-DaP20-case}
  \id-\gamma\geq \eta \text{ with a certain } \eta \in\Kinf,
\end{equation}
then $\gamma = \id - (\id-\gamma) \leq \id-\eta = (\id+\rho)^{-1}$ for a suitable $\rho\in\Kinf$ (this is easy to see, and was stated in \cite[Lem.~8]{MKG20}). In this case, we obtain that $\sigma$ is a path of strict decay as in Definition~\ref{def_omega_path}, and the corresponding ISS Lyapunov function constructed by our small-gain theorem is $V(s) = \sup_{i\in\N}V_i(s_i)$ for any $s = (s_i)_{i\in\N}\in\ell_\infty(\N,(n_i))$, which is precisely the Lyapunov function proposed in \cite[Thm.~1]{DaP19}. Thus, we obtain the main result of \cite{DaP19} as a special case of ours. Nevertheless, if the condition \eqref{eq:Paths-of-strict-decay-DaP20-case} does not hold, then $\Gamma$ decays at the points of the form $\sigma(r)$, but not in a uniform way as we require in this work. In this case, the result in \cite{DaP19} is applicable, while our theorem is not.%

%%%%%%%%%%%%%%%%

\textbf{Comparison to \cite{dashkovskiy2019stability}.} The setting of \cite{dashkovskiy2019stability} is very similar to that of our paper, in particular, the Lyapunov gains $\gamma_{ij}$ may be distinct nonlinear functions, and the number of neighbors is also taken to be finite. The small-gain theorem \cite[Thm.~5.1]{dashkovskiy2019stability} is shown under the requirement that there is a linear path of strict decay for the gain operator $\Gamma$, with a linear $\rho$ (see Definition~\ref{def_omega_path}), whereas we do not require linearity, which makes our result truly nonlinear.%

In \cite{dashkovskiy2019stability}, it was proposed to use the strong transitive closure $Q$ to construct a path of strict decay, and a robust SGC was introduced to describe under which conditions a path of strict decay exists. In \cite[Lem.~4.5]{dashkovskiy2019stability}, it was shown that there is $\sigma:\R_+ \to\ell_\infty^+$, satisfying properties (i)--(iii) of Definition~\ref{def_omega_path} with a linear $\rho$, provided that $\Gamma$ satisfies the robust strong SGC with linear $\rho,\omega$, and the discrete-time system \eqref{eq:Gammaotimes-discrete-time-nonlinear} induced by $\Gamma$ is UGAS. 

We improve these results in several directions. We introduce the concept of max-robust small-gain condition, which is weaker than the robust small-gain condition, and show that UGAS of the discrete-time system induced by the gain operator implies the max-robust small-gain condition (Proposition~\ref{prop_maxrobustsgc_char}). Finally, we show in Theorem~\ref{thm_omegapath}, that under certain regularity assumptions, UGAS of a system induced by a scaled gain operator is only needed to ensure the existence of a path of strict decay, which we also construct explicitly.

%%%%%%%%%%%%%%%

\textbf{Finite networks.} Finally, let us compare our results to the available results for finite networks. A nonlinear Lyapunov-based small-gain theorem for finite networks of ODE systems has been proposed in \cite{Dashkovskiy.2010}. As in the first part of this section, we require a bit more uniformity in the definition of a path of strict decay, and thus we cannot fully recover the small-gain result in \cite{Dashkovskiy.2010}. Also, since we formulate our result in a maximum formulation, the gain operator is a so-called max-preserving operator \cite{KaJ11,Rue17}, and it is known that the SGC (not even the strong one) is equivalent to UGAS of an induced discrete-time system and to the existence of a path of strict decay as defined in \cite{Dashkovskiy.2010}, see \cite[Thm.~6.4]{Rue10}.%

On the other hand, we know that for finite networks the strong SGC is equivalent to the robust strong SGC \cite[Prop.~14]{MKG20}, and implies UGAS of the discrete-time system \cite[Thm.~4.6]{Rue10}. Thus, for finite networks, our Theorem~\ref{thm_omegapath} states that the strong SGC for $\Gamma$ implies the existence of a path of strict decay with a uniform decay rate, characterized by $\rho$, which is not far from the sharp result for finite networks.%

\section{Conclusions and perspectives}

We have proved a fully nonlinear Lyapunov-based small-gain theorem for ISS of infinite networks. In our result, we use the maximum formulation and the implication form to describe the ISS property of the subsystems. A crucial assumption in our result is the existence of a (nonlinear) path of strict decay for the gain operator, which acts on the positive cone in the sequence space $\ell_{\infty}$. We have proved that such a path exists if the discrete-time system induced by a scaled gain operator is UGAS and, additionally, the interconnection gains satisfy a uniform local Lipschitz condition. While the second assumption is not hard to check in a concrete example, the first one is more delicate. We have characterized a weaker property in terms of the so-called max-robust small-gain condition, which is again tractable. The difference between UGAS and this weaker property is that trajectories may only converge to zero componentwise and not necessarily uniformly, i.e.~with respect to the $\ell_{\infty}$-norm. Different types of sufficient conditions for the uniform convergence were also provided and their tractability was shown in an example.%

One open question for future research is whether there exist characterizations of UGAS (of the gain-operator-induced system) that can be checked by a condition formulated in terms of single gains instead of their compositions. Another open problem is to extend the theory to different formulations of the ISS property for the subsystems, which lead to other types of gain operators.%

\section*{Appendix: Technical lemmas}

We omit the proofs of the following two simple lemmas, the first of which follows from \cite[Lem.~3.4]{Kha02}.%

\begin{lemma}\label{lem_techlem1}
Let $s_1,s_2:[0,T] \rightarrow \R_+$ be differentiable functions such that $s_1(0) = s_2(0)$, $\dot{s}_1(t) \leq -\alpha(s_1(t))$ and $\dot{s}_2(t) = -\alpha(s_2(t))$ for all $t \in [0,T]$ and some locally Lipschitz $\alpha \in \PC$. Then $s_1(t) \leq s_2(t)$ for all $t \in [0,T]$.
\end{lemma}

\begin{lemma}\label{lem_incrmax}
Let $f_i$, $i=1,\ldots,n$, be finitely many strictly increasing functions $f_i:[a,b] \rightarrow \R_+$ with $a,b \in \R$, $a < b$. Further assume that $|f_i(r_1) - f_i(r_2)| \geq l_i|r_1 - r_2|$ for all $r_1,r_2 \in [a,b]$ and $i \in \N$, where $l_i > 0$. Put $f(r) := \max_{i=1,\ldots,n}f_i(r)$, $f:[a,b] \rightarrow \R_+$. Then, with $l := \min\{l_1,\ldots,l_n\}$, we have%
\begin{equation*}
  |f(r_1) - f(r_2)| \geq l |r_1 - r_2| \mbox{\quad for all\ } r_1,r_2 \in [a,b].%
\end{equation*}
\end{lemma}

The following lemma can be found in \cite[Fact A.1.24]{downarowicz2011entropy}.%

\begin{lemma}\label{lem_supinf}
Let $(f_k)_{k\in\N}$ be a decreasing sequence of upper semicontinuous functions $f_k:X \rightarrow \R$ on a compact metric space $X$. Then%
\begin{equation*}
  \inf_{k\in\N}\sup_{x\in X}f_k(x) = \sup_{x\in X}\inf_{k\in\N}f_k(x).%
\end{equation*}
\end{lemma}

\bibliographystyle{IEEEtran}
\bibliography{references}

\begin{IEEEbiography}[{\includegraphics[width=1in,height=1.3in,clip,keepaspectratio]{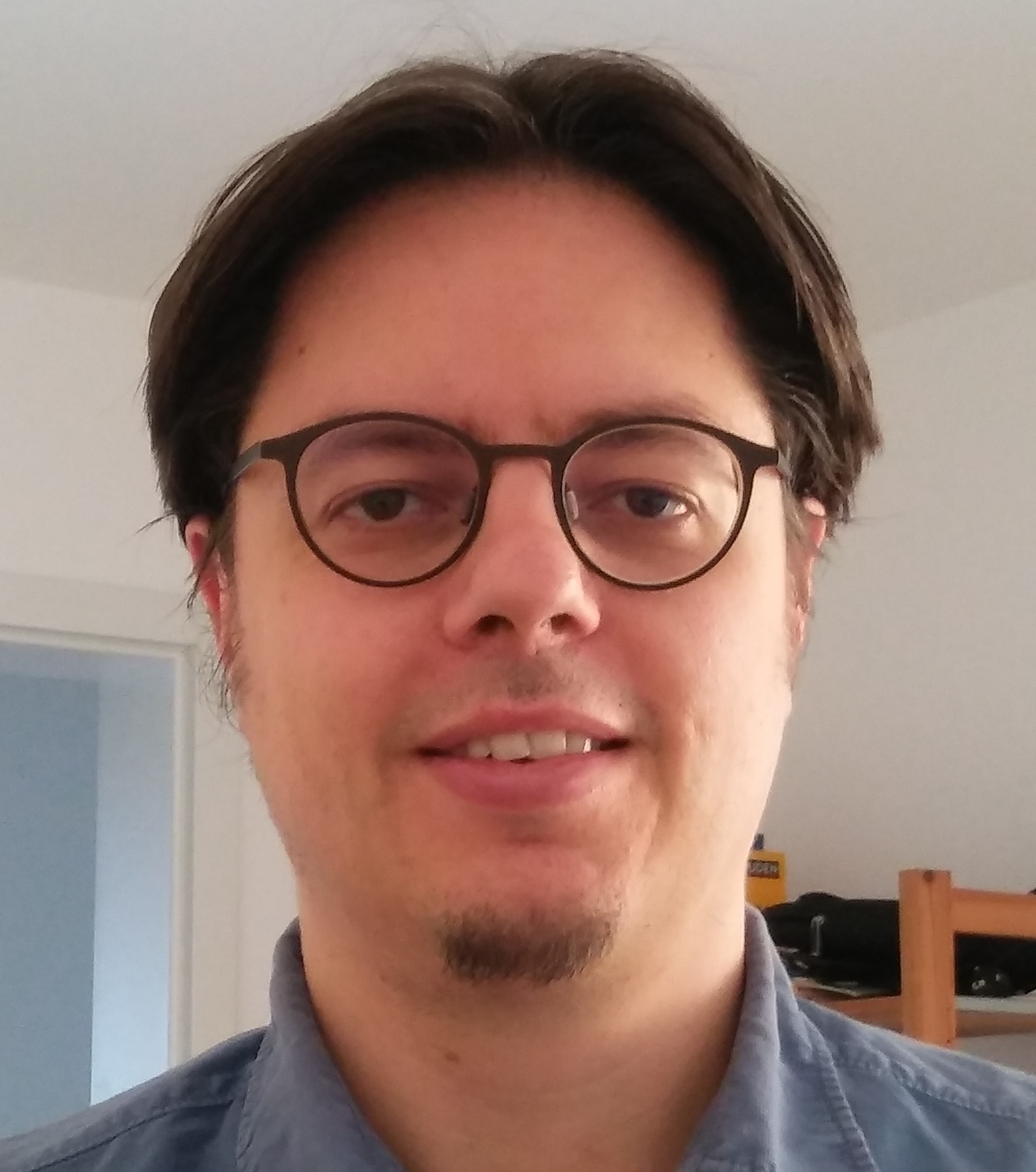}}]{Christoph Kawan} received the diploma and doctoral degree in mathematics (both under supervision of Prof. Fritz Colonius) 
from the University of Augsburg, Germany, in 2006 and 2009, respectively. He spent four months as a research scholar at the State University of Campinas, Brazil, in 2011, and nine months at the Courant Institute of Mathematical Sciences at New York University, in 2014. He is the author of the book `Invariance Entropy for Deterministic Control Systems - An Introduction' (Lecture Notes in Mathematics 2089. Berlin: Springer, 2013). His research interests include networked and information-based control and large-scale networks.%
\end{IEEEbiography}

\begin{IEEEbiography}[{\includegraphics[width=1in,height=1.25in,clip,keepaspectratio]{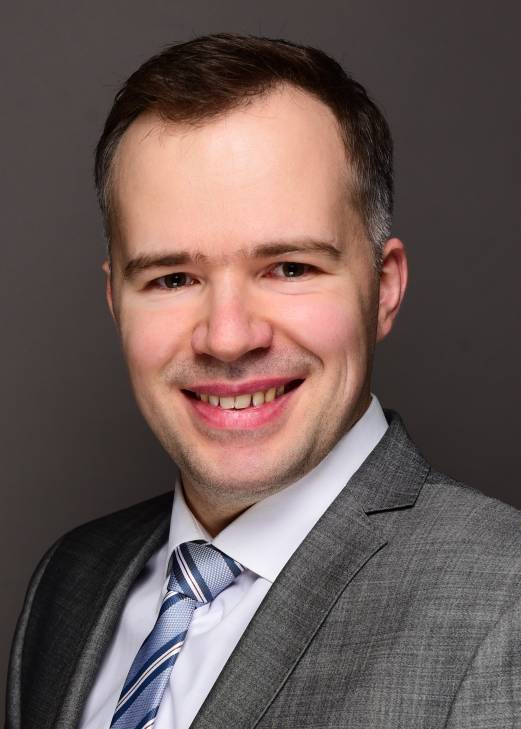}}]{Andrii Mironchenko}
was born in 1986. He received the M.Sc. degree in applied mathematics from the Odesa I.I. Mechnikov National University, Odesa, Ukraine, in 2008, and the Ph.D. degree in mathematics from the University of Bremen, Bremen, Germany in 2012. 

He has held a research position with the University of W\"urzburg, W\"urzburg, Germany and was a Postdoctoral Fellow of Japan Society for Promotion of Science (JSPS) with the Kyushu Institute of Technology, Fukuoka Prefecture, Japan (2013--2014). In 2014 he joined the Faculty of Mathematics and Computer Science, University of Passau, Passau, Germany.%

His research interests include infinite-dimensional systems, stability theory, hybrid systems and applications of control theory to biological systems and distributed control. 
\end{IEEEbiography}

\begin{IEEEbiography}[{\includegraphics[width=1in,height=1.25in,clip,keepaspectratio]{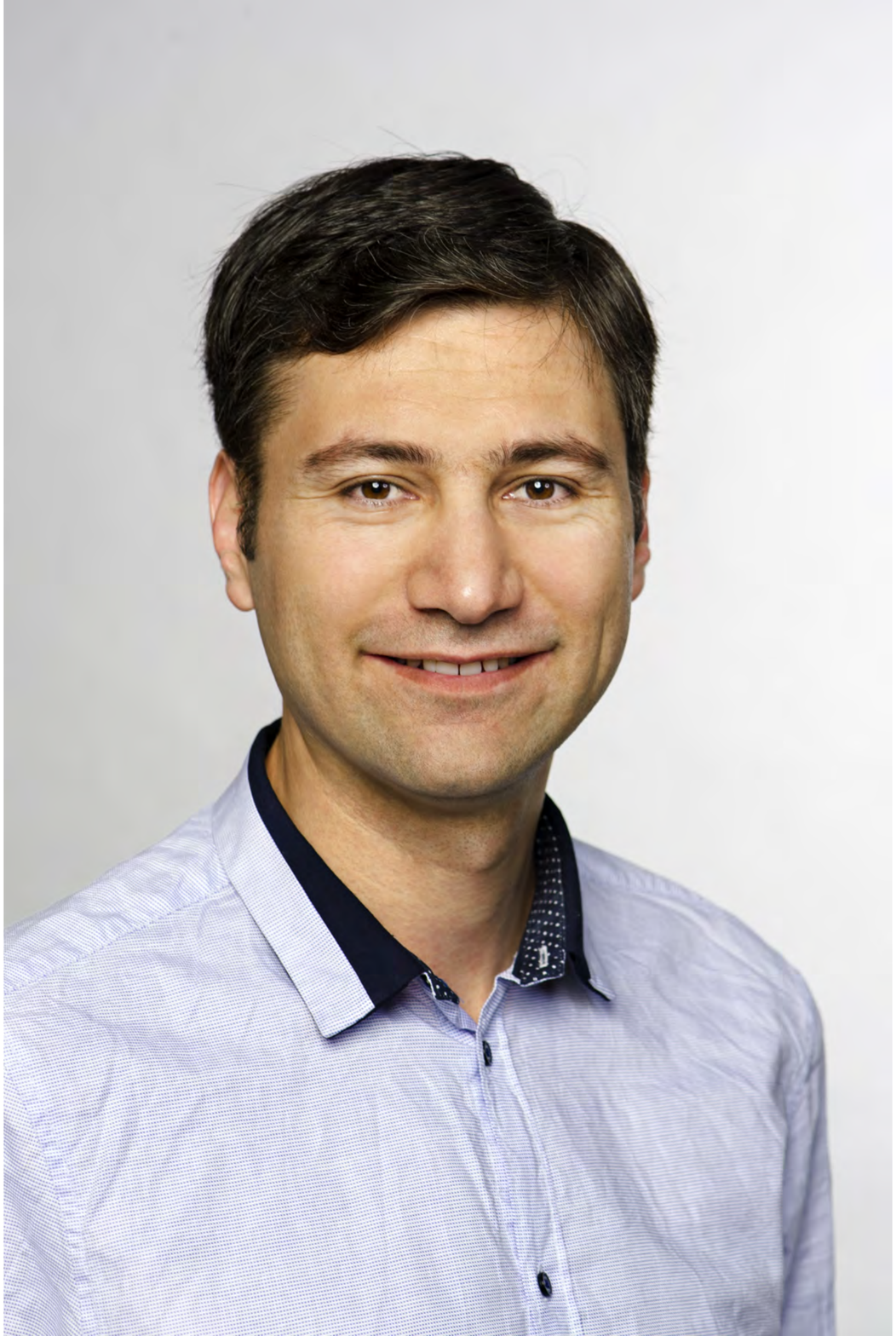}}]{Majid Zamani}
(M'12--SM'16)
is an Assistant Professor in the Computer Science Department at the University of Colorado Boulder, USA. He received a B.Sc. degree in Electrical Engineering in 2005 from Isfahan University of Technology, Iran, an M.Sc. degree in Electrical Engineering in 2007 from Sharif University of Technology, Iran, an MA degree in Mathematics and a Ph.D. degree in Electrical Engineering both in 2012 from University of California, Los Angeles, USA. Between September 2012 and December 2013 he was a postdoctoral researcher at the Delft Centre for Systems and Control, Delft University of Technology, Netherlands. From May 2014 to January 2019 he was an Assistant Professor in the Department of Electrical and Computer Engineering at the Technical University of Munich, Germany. From December 2013 to April 2014 he was an Assistant Professor in the Design Engineering Department, Delft University of Technology, Netherlands. He received an ERC starting grant award from the European Research Council in 2018.%

His research interests include verification and control of hybrid systems, embedded control software synthesis, networked control systems, and incremental properties of nonlinear control systems.
\end{IEEEbiography}

\end{document}